\documentclass[11pt]{article}
\usepackage{amsfonts}
\usepackage{euscript,amstext,a4}
\usepackage{mathrsfs}
\usepackage{indentfirst}
\usepackage{enumerate}
\usepackage{amssymb}
\usepackage{amsmath}
\usepackage{hyperref}
\usepackage{amscd}
\numberwithin{equation}{section}
\newtheorem{thm}{Theorem}[section]
\numberwithin{thm}{section}

\newtheorem{prop}[thm]{Proposition}
\newtheorem{rem}[thm]{Remark}
\newtheorem{lem}[thm]{Lemma}
\newtheorem{cor}[thm]{Corollary}

\DeclareMathOperator{\ind}{Ind}
\DeclareMathOperator{\sym}{Sym}

\title{On the Witten Rigidity Theorem for String$^c$ Manifolds}
\author{Jianqing YU\footnote{Chern Institute of Mathematics \& LPMC, Nankai
University, Tianjin 300071, P. R. China. (jianqingyu@gmail.com)} \ \ and\  \
Bo LIU\footnote{Chern Institute of Mathematics \& LPMC, Nankai University, Tianjin
300071, P. R. China. (boliumath@mail.nankai.edu.cn)}
}
\date{}

\begin{document}

\maketitle

\begin{abstract}
We establish the family rigidity and vanishing theorems on the equivariant $K$-theory level
for the Witten type operators on String$^c$ manifolds introduced by Chen-Han-Zhang \cite{chenhanzhang}.
\end{abstract}

\section{Introduction}
In \cite{MR970288}, Witten derived a series of elliptic operators on the free loop space $\mathcal{L}M$ of
a spin manifold $M$. In particular, the index of the formal signature operator on loop space turns out to be exactly
the elliptic genus constructed by Landweber-Stong \cite{MR948178} and Ochanine \cite{MR895567} in a topological way. Motivated by
physics, Witten proposed that these elliptic operators should be rigid with respect to the circle action.

This claim of Witten was first proved by Taubes \cite{MR998662} and Bott-Taubes \cite{MR954493}.
See also \cite{MR981372} and \cite{MR1048541} for other interesting cases.
By the modular invariance property, Liu (\cite{MR1331972,MR1396769}) presented a simple
and unified proof of the above result as well as various further generalizations.
In particular, several new vanishing theorems were established in \cite{MR1331972,MR1396769}.

In a recent paper \cite{chenhanzhang}, Chen, Han and Zhang introduced a topological condition
which they called String$^c$ condition for even dimensional Spin$^c$ manifolds. Under this
String$^c$ condition, they constructed a Witten type genus which is the index of a Witten
type operator, a linear combination of twisted Spin$^c$ Dirac operators. Furthermore, by applying Liu's method in
\cite{MR1331972,MR1396769}, Chen-Han-Zhang established the rigidity and vanishing theorems for
this Witten type operator under relevant anomaly cancelation condition (cf. \cite[Theorem 3.2]{chenhanzhang}).

In many situations in geometry, it is rather natural and necessary to generalize the
rigidity and vanishing theorems to the family case.
On the equivariant Chern character level, Liu and Ma (\cite{MR1756105,MR1969037})
established several family rigidity and vanishing theorems.
In \cite{MR1870666,MR2016198}, inspired by \cite{MR998662}, Liu, Ma
and Zhang established the corresponding family
rigidity and vanishing theorems on the equivariant $K$-theory level.
As explained in \cite{MR1870666,MR2016198}, taking the Chern character might
kill some torsion elements involved in the index bundle.  Therefore, the $K$-theory level
rigidity and vanishing properties are more subtle than those on the
Chern character level.

The purpose of this paper is to establish the family rigidity and vanishing theorems on the equivariant $K$-theory level
for the Witten type operators introduced by Chen-Han-Zhang \cite{chenhanzhang}.
In fact, our main results in Theorem \ref{main} may be regarded as
an analogue of \cite[Theorem 2.1]{MR1870666}
and \cite[Theorems 2.1, 2.2]{MR2016198}. In particular, if the base
manifold is a point, from our family rigidity theorem, one deduces Chen-Han-Zhang's
theorem \cite[Theorem 3.2(i)]{chenhanzhang}.
Both the statement and the proof of Theorem \ref{main} are
inspired by those of \cite[Theorem 2.1]{MR1870666} and
\cite[Theorems 2.1, 2.2]{MR2016198}, which essentially
depend on the techniques developed by Taubes \cite{MR998662} and Bismut-Lebeau \cite{MR1188532}.

This paper is organized as follows. In Section \ref{sec2}, we state and prove our main results, Theorem \ref{main},
the rigidity and vanishing theorems for the family Witten type operators introduced by Chen-Han-Zhang \cite{chenhanzhang}.
Section \ref{sec3} is devoted to the proofs of two intermediate results, Theorems \ref{eg31} and \ref{eg29}, which are used in
the proof of Theorem \ref{main}.

\section{Rigidity and vanishing theorems in $K$-theory}\label{sec2}

In this section, we establish the main results of this paper, the rigidity and vanishing theorems
on the equivariant $K$-theory level for a family of Spin$^c$ manifolds. Such theorems hold under some
anomaly cancelation assumption which is inspired by Chen-Han-Zhang's String$^c$ condition \cite{chenhanzhang}.
For the particular case when the base manifold is a point, our results imply Chen-Han-Zhang's
theorem \cite[Theorem 3.2(i)]{chenhanzhang}.

This section is organized as follows. In Section \ref{sec2.1}, we reformulate
a $K$-theory version of the equivariant family index theorem which is proved in
\cite[Theorem 1.2]{MR2016198} and \cite[Theorem 1.1]{MR1870666}. In Section \ref{sec2.2},
we state our main results, the rigidity and vanishing theorems
on the equivariant $K$-theory level for a family of Spin$^c$ manifolds.
In Section \ref{sec2.3}, we state two intermediate results
on the relations between the family indices on the fixed point set, which
will be used to prove our main results stated in Section \ref{sec2.1}.
In Section \ref{sec2.4}, we prove the family rigidity and vanishing theorems.

\subsection{A $K$-theory version of the equivariant family index theorem}\label{sec2.1}

Let $M$, $B$ be two compact manifolds, and $\pi:M\rightarrow B$ a smooth fibration with compact fiber $X$ such that
$\dim X=2l$. Let $TX$ denote the relative tangent bundle carrying a Riemannian metric $g^{TX}$. We assume that
$TX$ is oriented. Let $(W,h^W)$ be a complex Hermitian vector bundle over $M$.

Let $(V,g^{V})$ (resp. $(V',g^{V'})$) be a $2p$ (resp. $2p'$) dimensional oriented
real Euclidean vector bundle over $M$. Let $(L,h^L)$ be a complex Hermitian line bundle over $M$ with the property that the
vector bundle $U=TX\oplus V\oplus V'$ satisfies $\omega_2(U)=c_1(L)\mod (2)$,
where $\omega_2$ denotes the second Stiefel-Whitney
class, and $c_1$ denotes the first Chern class. Then the  vector
bundle $U$ has a Spin$^c$-structure.
Let $S(U,L)$ be the fundamental complex spinor bundle for $(U,L)$ (cf. \cite[Appendix D]{MR1031992}).

Assume that there is a fiberwise $S^1$ action on $M$ which
lifts to $V$, $V'$, $L$ and $W$, and assume the metrics
$g^{TX}$, $g^{V}$, $g^{V'}$, $h^L$ and $h^W$ are $S^1$-invariant. Also assume that the $S^1$ actions
on $TX$, $V$, $V'$, $L$ lift to $S(U,L)$.

Let $\nabla^{TX}$ be the Levi-Civita connection on $(TX,g^{TX})$ along the fiber $X$.
Let $\nabla^{V}$ (resp. $\nabla^{V'}$) be an $S^1$-invariant Euclidean connection on $(V,g^{V})$ (resp. $(V',g^{V'})$).
Let $\nabla^{L}$ (resp. $\nabla^{W}$) be an $S^1$-invariant Hermitian connection on $(L,h^L)$ (resp. $(W,h^W)$).

The Clifford algebra bundle $C(TX)$ is the bundle of Clifford algebras over $X$ whose fibre at $x\in X$ is the
Clifford algebra $C(T_xX)$ (cf. \cite{MR1031992}). Let $C(V)$ (resp. $C(V')$) be the Clifford algebra bundle of
$(V,g^V)$ (resp. $(V',g^{V'})$).

Let $\{e_i\}_{i=1}^{2l}$ (resp. $\{f_j\}_{j=1}^{2p}$) be an
oriented orthonormal basis of $(TX, g^{TX})$ (resp. $(V, g^{V})$).
We denote by $c(\cdot)$ the Clifford action of $C(TX)$, $C(V)$ and $C(V')$ on $S(U,L)$.
Let $\tau$ be the involution of $S(U,L)$ given by
\begin{equation}\label{s5}
\tau=(\sqrt{-1})^{l+p}c(e_1)\cdots c(e_{2l})c(f_1)\cdots c(f_{2p})\ .
\end{equation}
In the rest of the paper, we will say $\tau$ the involution determined by $TX\oplus V$.
We decompose $S(U,L)=S_+(U,L)\oplus S_-(U,L)$ corresponding to $\tau$
such that $\tau\big|_{S_{\pm}(U,L)}=\pm1$. Let $\nabla^{S(U,L)}$ be the Hermitian connection on $S(U,L)$ induced by
$\nabla^{TX}$, $\nabla^{V}$, $\nabla^{V'}$ and $\nabla^L$ (cf. \cite[Appendix D]{MR1031992}).
Then $\nabla^{S(U,L)}$ preserves the $\mathbb{Z}_2$-grading of $S(U,L)$ induced by \eqref{s5}.
Let $\nabla^{S(U,L)\otimes W}$ be the Hermitian connection on
$S(U,L)\otimes W$ obtained from the
tensor product of $\nabla^{S(U,L)}$ and $\nabla^W$.
Let $D^{X}\otimes W$ be the family twisted Spin$^c$-Dirac operator on the fiber $X$
defined by
\begin{equation}\label{s22}
D^{X}\otimes W=\sum_{i=1}^{2l}c(e_i)\nabla^{S(U,L)\otimes W}_{e_i}\ .
\end{equation}
By \cite[Proposition 1.1]{MR1756105}, the index bundle $\ind_\tau(D^X\otimes W)$ over $B$ is well-defined
in the equivariant $K$-group $K_{S^1}(B)$. Using the same notations as in \cite[(1.4)-(1.7)]{MR2016198}, we
write, as an identification of virtual $S^1$-bundles,
\begin{equation}
\ind_\tau(D^X\otimes W)=\bigoplus_{n\in \mathbb{Z}}\ind_\tau(D^X\otimes W,n)\otimes [n],
\end{equation}
where by $[n]$ ($n\in\mathbb{Z}$) we mean the one dimensional complex vector space on which $S^1$ acts
as multiplication by $g^n$ for a generator $g\in S^1$.

Let $F=\{F_\alpha\}$ be the fixed point set of the circle action on $M$. Then $\pi:F_\alpha\rightarrow B$ (resp. $\pi:F\rightarrow B$)
is a smooth fibration with fiber $Y_\alpha$ (resp. $Y$). Let $\widetilde{\pi}:N\rightarrow F$ denote the
normal bundle to $F$ in $M$. Then $N=TX/TY$. We identify $N$ as the orthogonal complement of $TY$ in $TX|_{F}$.
Then $TX|_{F}$ admits a $S^1$-equivariant decomposition (cf. \cite[(1.8)]{MR2016198})
\begin{equation}\label{s19}
TX|_{F}=\bigoplus_{v\neq 0}N_v\oplus TY\,,
\end{equation}
where $N_v$ is a complex vector bundle such that $g\in S^1$ acts on it by $g^v$ with $v\in \mathbb{Z}\backslash \{0\}$.
Clearly, $N=\oplus_{v\neq 0}N_v$. We will regard $N$ as a complex vector bundle and write
$N_\mathbb{R}$ for the underlying real vector bundle of $N$. For $v\neq 0$, let $N_{v,\mathbb{R}}$ denote
the underlying real vector bundle of $N_v$.

Similarly, let (cf. \cite[(1.9) and (1.46)]{MR2016198})
\begin{equation}\label{s20}
V\big|_{F}=\bigoplus_{v\neq 0}V_{v}\oplus V_{0}^{\mathbb{R}}\ ,
\quad {V'}|_{F}=\bigoplus_{v\neq 0}V'_v\oplus V_0^{\prime\mathbb{R}}\ ,
\quad W|_{F}=\bigoplus_{v}W_v\ ,
\end{equation}
be the $S^1$-equivariant decompositions of the restrictions of $V$, $V'$ and $W$ over $F$ respectively, where
$V_v$, $V'_v$ and $W_v$ $(v\in \mathbb{Z})$ are  complex vector bundles over $F$
on which $g\in S^1$ acts by $g^v$, and $V_0^{\mathbb{R}}$ (resp. $V_0^{\prime\mathbb{R}}$)
is the real subbundle of $V$ (resp. $V'$) such that
$S^1$ acts as identity. For $v\neq 0$, let $V_{v,\mathbb{R}}$ (resp. $V'_{v,\mathbb{R}}$) denote
the underlying real vector bundle of $V_v$ (resp. $V'_v$). Denote by
$2p_0=\dim V_0^{\mathbb{R}}$ and $2l_0=\dim Y$.

Let us write (compare with \cite[(1.47)]{MR2016198})
\begin{align}
L_{F}=L\otimes\Bigl(\,\bigotimes_{v\neq0}\det
N_v\otimes\bigotimes_{v\neq0}\det V_v\otimes\bigotimes_{v\neq0}\det V'_v\,\Bigr)^{-1}.
\end{align}
Then $TY\oplus V_0^{\mathbb{R}}\oplus V_0^{\prime\mathbb{R}}$ has a Spin$^c$-structure.
Let $S(TY\oplus V_0^{\mathbb{R}}\oplus V_0^{\prime\mathbb{R}},L_F)$ be the fundamental spinor bundle for
$(TY\oplus V_0^{\mathbb{R}}\oplus V_0^{\prime\mathbb{R}},L_F)$. Let $R$ be a Hermitian complex vector
bundle equipped with a Hermitian connection over $F$.
We will denote by $D^{Y}\otimes R$ the family (twisted) Spin$^c$
Dirac operator on $S(TY\oplus V_0^{\mathbb{R}}\oplus V_0^{\prime\mathbb{R}},L_F)\otimes R$ defined as in \eqref{s22}
and by $D^{Y_\alpha}\otimes R$ its restriction to $Y_\alpha$.

Recall that $N_{v,\mathbb{R}}$ and $V_{v,\mathbb{R}}$ are canonically oriented by
their complex structures. The decompositions (\ref{s19}),
(\ref{s20}) induce the orientations of $TY$ and $V_{0}^{\mathbb{R}}$ respectively. Let
$\{e_i\}_{i=1}^{2l_0}$, $\{f_j\}_{j=1}^{2p_0}$ be the corresponding
oriented orthonormal basis of $(TY,g^{TY})$ and
$(V_0^{\mathbb{R}},g^{V_0^{\mathbb{R}}})$. The involution of $S(TY\oplus V_0^{\mathbb{R}}\oplus V_0^{\prime\mathbb{R}},L_F)$ is
canonically associated to that of $S(U,L)$, which we still denote by $\tau$, is given by
\begin{equation}
\tau=(\sqrt{-1})^{l_0+p_0}c(e_1)\cdots c(e_{2l_0})c(f_1)\cdots c(f_{2p_0})\ .
\end{equation}
Let $S(TY\oplus V_0^{\mathbb{R}}\oplus V_0^{\prime\mathbb{R}},L_F)
=S_{+}(TY\oplus V_0^{\mathbb{R}}\oplus V_0^{\prime\mathbb{R}},L_F)
\oplus S_{-}(TY\oplus V_0^{\mathbb{R}}\oplus V_0^{\prime\mathbb{R}},L_F)$ be the $\mathbb{Z}_2$-grading of
$S(TY\oplus V_0^{\mathbb{R}}\oplus V_0^{\prime\mathbb{R}},L_F)$ induced by $\tau$.

Let $C(N_\mathbb{R})$ (resp. $C(V_{v,\mathbb{R}})$) be the
Clifford algebra bundle of $(N_\mathbb{R},g^{TX}|_{N_\mathbb{R}})$ (resp. $(V_{v,\mathbb{R}},g^{V}|_{V_{v,\mathbb{R}}})$).
By \cite[(1.10)]{MR2016198},  $\Lambda(\overline{N}^*)$ is a $C(N_\mathbb{R})$-Clifford module
with the involution
$\tau^N|_{\Lambda^{{\rm even}/{\rm odd}}(\overline{N}^*)}=\pm 1$.
Similarly as in \cite[(1.10)]{MR2016198}, we can define the Clifford action of $C(V_{v,\mathbb{R}})$
on $\Lambda(\overline{V_v}^*)$. Then $\Lambda(\overline{V_v}^*)$ is a $C(V_{v,\mathbb{R}})$-Clifford module
with the involution $\tau^V_v\big|_{\Lambda^{{\rm even}/{\rm odd}}(\overline{V_v}^*)}=\pm 1$.

By restricting to $F$, one has the isomorphism of
$\mathbb{Z}_2$-graded 
$C(TX)$-Clifford modules over $F$ as follows (compare with \cite[(1.49)]{MR2016198}),
\begin{equation}
\begin{split}\label{s21}
\bigl(S(U,L),\tau\bigr)\big|_{F}\simeq
&\ \Bigl(S(TY\oplus V_0^{\mathbb{R}}\oplus V_0^{\prime \mathbb{R}},L_F),\tau\Bigr)
\widehat{\otimes}\,\bigl(\Lambda\overline{N}^*,\tau^N\bigr)
\\
&\quad \widehat{\otimes}\,\widehat{\bigotimes_{v\neq0}}\,\bigl(\Lambda \overline{V_v}^*,\tau_v^V\bigr)\widehat{\otimes}\,\widehat{\bigotimes_{v\neq0}}\,\bigl(\Lambda \overline{V'_v}^*,{\rm id}\bigr),
\end{split}
\end{equation}
where ${\rm id}$ denotes the trivial involution and $\widehat{\otimes}$ denotes the $\mathbb{Z}_2$-graded tensor product
(cf. \cite[pp. 11]{MR1031992}). Furthermore, the
isomorphism  \eqref{s21} gives the identifications of the canonical connections on
the bundles (compare with \cite[(1.13)]{MR2016198}).

Let $S^1$ act on $L|_{F}$ by sending $g\in S^1$ to $g^{l_c}$ ($l_c\in \mathbb{Z}$) on $F$.
Then $l_c$ is locally constant on $F$. Following \cite[(1.50)]{MR2016198}, we define the following elements in $K(F)[[q^{\frac{1}{2}}]]$,
\begin{align}
R(q)&=q^{\frac{1}{2}\left(\sum_v|v|\dim N_v-\sum_v v\dim V_{v}-\sum_v v\dim V'_{v}+l_c\right)}
\bigotimes_{v>0}\left({\rm
Sym}_{q^v}(N_v)\otimes\det N_v\right)
\notag\\
&\quad\otimes\bigotimes_{v<0}{\rm
Sym}_{q^{-v}}(\overline{N}_v)\otimes\bigotimes_{v\neq 0}\Lambda_{-q^v}(V_{v})
\otimes\bigotimes_{v\neq 0}\Lambda_{q^v}(V'_{v})\otimes\Bigl(\,\sum_vq^vW_v\Bigr)
\notag\\
&=\sum_{n} R_{n}q^n\ ,
\\
R'(q)&=q^{\frac{1}{2}\left(-\sum_v|v|\dim N_v-\sum_v v\dim V_{v}-\sum_v v\dim V'_{v}+l_c\right)}
\bigotimes_{v>0}{\rm
Sym}_{q^{-v}}(\overline{N}_v)
\notag\\
&\quad \otimes\bigotimes_{v<0}\left({\rm
Sym}_{q^v}(N_v)\otimes\det N_v\right)\otimes\bigotimes_{v\neq 0}\Lambda_{-q^v}(V_{v})
\notag\\
&\quad\otimes\bigotimes_{v\neq 0}\Lambda_{q^v}(V'_{v})\otimes\Bigl(\,\sum_vq^vW_v\Bigr)
=\sum_{n} R'_{n}q^n.
\end{align}

As explained in \cite[pp. 139]{MR2016198}, since $TX\oplus V\oplus V'\oplus L$ is spin,
one gets
\begin{equation}
\sum_vv\dim N_v+\sum_v v\dim V_{v}+\sum_v v\dim V'_{v}+l_c \equiv 0 \mod (2).
\end{equation}
Therefore, $R_{n}(q)$, $R'_{n}(q)\in K(F)[[q]]$.

The following theorem was essentially proved in \cite[Theorem 1.2]{MR2016198}.

\begin{thm}\label{local}
For $n \in\mathbb{Z}$, the following identity holds in $K(B)$,
\begin{equation}
\begin{split}\label{s43}
\ind_{\tau}\bigl(D^X\otimes W,n\bigr)
&=\sum_{\alpha}(-1)^{\sum_{0<v}\dim N_v}\ind_{\tau}\bigl(D^{Y_\alpha}\otimes R_{n}\bigr)
\\
&=\sum_{\alpha}(-1)^{\sum_{v<0}{\rm
dim}\,N_v}\ind_{\tau}\bigl(D^{Y_\alpha}\otimes R'_{n}\bigr)\,.
\end{split}
\end{equation}
\end{thm}

\subsection{Family rigidity and vanishing theorems}\label{sec2.2}

Let $\pi:M\rightarrow B$ be a fibration of compact manifolds with fiber $X$ and $\dim X=2l$.
We assume that $S^1$ acts fiberwise on $M$ and $TX$ has an $S^1$-invariant Spin$^c$ structure.
Let $K_X$ be the $S^1$-equivariant complex line bundle over $M$ which is induced by the $S^1$-invariant
Spin$^c$ structure of $TX$. Let $S(TX,K_X)$ be the complex spinor bundle of $(TX,K_X)$ (cf. \cite[Appendix D]{MR1031992}).

Let $V$ be an even dimensional real vector bundle over $M$. We assume that $V$ has
an $S^1$-invariant spin structure. Let $S(V)=S^+(V)\oplus S^{-}(V)$ be the spinor bundle of $V$.
Let $W$ be an $S^1$-equivariant complex vector bundle
over $M$. Let $K_W=\det (W)$ be the determinant line bundle of $W$.

We define the following elements in $K(M)[[q^{1/2}]]$,
\begin{equation}\label{s3}
\begin{split}
R_1(V)&=\Bigl(S^+(V)+S^-(V)\Bigr)\otimes\bigotimes_{n=1}^{\infty}\Lambda_{q^n}(V)\ ,
\\
R_2(V)&=\Bigl(S^+(V)-S^-(V)\Bigr)\otimes\bigotimes_{n=1}^{\infty}\Lambda_{-q^n}(V)\ ,
\\
R_3(V)&=\bigotimes_{n=1}^{\infty}\Lambda_{-q^{n-1/2}}(V)\ ,\qquad
R_4(V)=\bigotimes_{n=1}^{\infty}\Lambda_{q^{n-1/2}}(V)\ ,
\end{split}
\end{equation}
\begin{equation*}
\begin{split}\label{s1}
Q_1(W)=&\bigotimes_{n=0}^\infty\Lambda_{q^n}(\overline{W})\otimes
\bigotimes_{n=1}^\infty\Lambda_{q^n}(W)\otimes
\bigotimes_{n=1}^\infty\Lambda_{-q^{n-1/2}}(\overline{W})
\\
&\ \otimes\bigotimes_{n=1}^\infty\Lambda_{-q^{n-1/2}}(W)
\otimes
\bigotimes_{n=1}^\infty\Lambda_{q^{n-1/2}}(\overline{W})\otimes
\bigotimes_{n=1}^\infty\Lambda_{q^{n-1/2}}(W)\ .
\end{split}
\end{equation*}

For $N\in \mathbb{Z}$, $N\geq 1$, let $y=e^{2\pi i/N}\in \mathbb{C}$. Let $G_y$ be the
multiplicative group generated by $y$. Following \cite{MR970288}, as in \cite[Section 2.1]{MR1870666},
we consider the fiberwise action $G_y$ on $W$ and $\overline{W}$ by sending $y\in G_y$ to $y$ on $W$ and
$y^{-1}$ on $\overline{W}$. Then $G_y$ acts naturally on $Q_1(W)$.

Let $H^*_{S^1}(M,\mathbb{Z})=H^*(M\times_{S^1}ES^1,\mathbb{Z})$ denote the $S^1$-equivariant cohomology group
of $M$, where $ES^1$ is the universal $S^1$-principal bundle over the classifying space $BS^1$ of $S^1$. So
$H^*_{S^1}(M,\mathbb{Z})$ is a module over $H^*(BS^1,\mathbb{Z})$ induced by the projection
$\overline{\pi}:M\times_{S^1}ES^1\rightarrow BS^1$.
Let $p_1(\cdot)_{S^1}$ and $\omega_2(\cdot)_{S^1}$ denote the first $S^1$-equivariant pontrjagin class and
the second $S^1$-equivariant Stiefel-Whitney class, respectively. As $V\times_{S^1}ES^1$ is spin over $M\times_{S^1}ES^1$, one knows that $\frac{1}{2}p_1(V)_{S^1}$ is well defined in $H^*_{S^1}(M,\mathbb{Z})$ (cf. \cite[pp. 456-457]{MR998662}). Recall that
\begin{equation}\label{eg5}
H^*(BS^1,\mathbb{Z})=\mathbb{Z}[[u]]
\end{equation}
with $u$ a generator of degree 2.

In the following, we denote by $D^X\otimes R$ the family twisted Spin$^c$ Dirac operator acting
fiberwise on $S(TX,K_X)\otimes R$.
Recall that if $\ind (D^X\otimes R,n)$ vanishes for $n\neq 0$, we say that $D^X\otimes R$ is rigid
on the equivariant $K$-theory level for the $S^1$ action.

Now we can state the main results of this paper as follows, which can
be thought of as an analogue of \cite[Theorem 2.1]{MR1870666}.

\begin{thm}\label{main}
Assume $w_2(W)_{S^1}=w_2(TX)_{S^1}$, $\frac{1}{2}p_1(V+3W-TX)_{S^1}=e\cdot\overline{\pi}^*u^2$ $(e\in\mathbb{Z})$
in $H^*_{S^1}(M,\mathbb{Z})$, and $c_1(W)=0 \mod (N)$. For $i=1,2,3,4$,
consider the family of $G_y\times S^1$-equivariant twisted Spin$^c$ Dirac operators
\begin{equation}
D^X\otimes(K_W\otimes K^{-1}_X)^{1/2}\otimes\bigotimes_{n=1}^{\infty}\sym_{q^n}(TX)\otimes R_i(V)\otimes Q_1(W)\ .
\end{equation}
\begin{enumerate}[{\rm (i)}]
\item If $e=0$, then these operators are rigid on the equivariant $K$-theory level for the $S^1$ action.
\item If $e<0$, then the index bundles of these operators are zero in $K_{G_y\times S^1}(B)$. In particular,
these index bundles are zero in $K_{G_y}(B)$.
\end{enumerate}
\end{thm}

\begin{rem}
As explained in \cite[Remark 2.1]{MR1870666}, $w_2(W)_{S^1}=w_2(TX)_{S^1}$ means that
$\frac{1}{2}p_1(3W-TX)_{S^1}$ is well defined and that $c_1(K_W\otimes K^{-1}_X)_{S^1}=0\mod(2)$.
By \cite[Corollary 1.2]{MR0461538}, the $S^1$ action on $M$ can be lifted to $(K_W\otimes K^{-1}_X)^{1/2}$
and is compatible with the $S^1$ action on $K_W\otimes K_X^{-1}$.
\end{rem}

Take $N=1$, i.e., we forget the $G_y$ action on $W$ and remove the corresponding assumption $c_1(W)=0 \mod (N)$.
Furthermore, take $W=K_X$ and $V=0$. Then
an interesting consequence of Theorem \ref{main} is the following family rigidity and vanishing property,
which may be thought of as an extension of \cite[Theorem 2.3]{MR2016198} to the Spin$^c$ case.
When the base manifold is a point, it turns out exactly to be Chen-Han-Zhang's theorem \cite[Theorem 3.2(i)]{chenhanzhang}.

\begin{cor}
Assume $\frac{1}{2}p_1(3K_X-TX)_{S^1}=e\cdot\overline{\pi}^*u^2$ $(e\in\mathbb{Z})$
in $H^*_{S^1}(M,\mathbb{Z})$. Consider the family of $S^1$-equivariant twisted Spin$^c$ Dirac operators
\begin{equation}\label{s4}
D^X\otimes\bigotimes_{n=1}^{\infty}\sym_{q^n}(TX)\otimes Q_1(K_X)\ .
\end{equation}
\begin{enumerate}[{\rm (i)}]
\item If $e=0$, then these operators are rigid on the equivariant $K$-theory level for the $S^1$ action.
\item If $e<0$, then the index bundles of these operators are zero in $K_{S^1}(B)$. In particular,
these index bundles are zero in $K(B)$.
\end{enumerate}
\end{cor}

\begin{rem}
The operators in \eqref{s4} are the Witten type operators introduced by Chen-Han-Zhang \cite{chenhanzhang}.
By taking $N=1$, $W=K_X$, $V=0$, and letting the base manifold $B$ be a point in \cite[Theorem 2.1]{MR1870666},
we get Chen-Han-Zhang's theorem \cite[Theorem 3.2(ii)]{chenhanzhang}.
It is rather natural to establish an analogue of
\cite[Theorem 2.1]{MR1870666}, which corresponds to Chen-Han-Zhang's
theorem \cite[Theorem 3.2(i)]{chenhanzhang}. That is
one of the motivation of Theorem \ref{main}.
\end{rem}

Actually, as in \cite{MR1870666,MR2016198}, our proof of Theorem \ref{main} works under the following slightly weaker hypothesis.
Let us first explain some notations.

For each $n>1$, consider $\mathbb{Z}_n\subset S^1$, the cyclic subgroup of order $n$. We have the $\mathbb{Z}_n$-equivariant
cohomology of $M$ defined by $H^*_{\mathbb{Z}_n}(M,\mathbb{Z})=H^*(M\times_{\mathbb{Z}_n}ES^1,\mathbb{Z})$,
and there is a natural \textquoteleft\textquoteleft forgetful\textquoteright\textquoteright\ map $\alpha(S^1,\mathbb{Z}_n):M\times_{\mathbb{Z}_n}ES^1\rightarrow M\times_{S^1}
ES^1$ which induces a pullback $\alpha(S^1,\mathbb{Z}_n)^*:H^*_{S^1}(M,\mathbb{Z})\rightarrow
H^*_{\mathbb{Z}_n}(M,\mathbb{Z})$. We
denote by $\alpha(S^1,1)$ the arrow which forgets the $S^1$ action.
Thus $\alpha(S^1,1)^*:H^*_{S^1}(M,\mathbb{Z})\rightarrow H^*(M,\mathbb{Z})$
is induced by the inclusion of $M$ into $M\times_{S^1}ES^1$ as a fiber over $BS^1$.

Finally, note that if $\mathbb{Z}_n$ acts trivially on a space $Y$, then there is a new arrow $t^*:H^*(Y,\mathbb{Z})
\rightarrow H^*_{\mathbb{Z}_n}(Y,\mathbb{Z})$ induced by the projection
$t:Y\times_{\mathbb{Z}_n}ES^1=Y\times B\mathbb{Z}_n\rightarrow Y$.

Let $\mathbb{Z}_{\infty}=S^1$. For each $1<n\leq+\infty$, let $i:M(n)\rightarrow M$ be the
inclusion of the fixed point set of $\mathbb{Z}_n\subset S^1$ in $M$, and so $i$ induces
$i_{S^1}:M(n)\times_{S^1}ES^1\rightarrow M\times_{S^1}ES^1$.

In the rest of this paper, we suppose that there exists some
integer $e\in\mathbb{Z}$ such that for $1<n\leq+\infty$,
\begin{equation}\label{hypothesis}
\begin{split}
\alpha(S^1,\mathbb{Z}_n)^*\circ i_{S^1}^*\Bigl(\,\frac{1}{2}p_1(V+3W-TX)_{S^1}-e\cdot\overline{\pi}^*u^2\Bigr)&
\\
=t^*\circ\alpha(S^1,1)^*\circ i_{S^1}^*\Bigl(\,\frac{1}{2}p_1(V+3W-TX)_{S^1}\Bigr)\ .&
\end{split}
\end{equation}
As indicated in \cite[Remark 2.4]{MR1870666}, the relation (\ref{hypothesis})
clearly follows from the hypothesis of Theorem \ref{main} by
pulling back and forgetting. Thus it is a weaker hypothesis.

We can now state a slightly more general version of Theorem \ref{main}.
\begin{thm}\label{theorem}
Under the hypothesis {\rm (\ref{hypothesis})}, we have
\begin{enumerate}[{\rm (i)}]
\item If $e=0$, then the index bundles of the twisted Spin$^c$ Dirac operators in
Theorem \ref{main} are rigid on the equivariant $K$-theory level for the $S^1$ action.

\item If $e<0$, then the index bundles of the twisted Spin$^c$ Dirac operators in
Theorem \ref{main} are zero as elements in $K_{G_y\times S^1}(B)$. In particular,
these index bundles are zero in $K_{G_y}(B)$.
\end{enumerate}
\end{thm}

The rest of this section is devoted to a proof of Theorem \ref{theorem}.

\subsection{Two recursive formulas}\label{sec2.3}

Let $F=\{F_\alpha\}$ be the fixed point set of the circle action. Then
$\pi:F\rightarrow B$ is a fibration with compact fibre denoted by $Y=\{Y_\alpha\}$.

As in \cite[(2.5)]{MR1870666},
we may and we will assume that
\begin{equation}\label{eg8}
\begin{split}
&TX|_{F}=TY\oplus\bigoplus_{v>0}N_v,\\
&TX|_{F}\otimes_{\mathbb{R}}\mathbb{C}=TY\otimes_{\mathbb{R}}\mathbb{C}
\oplus\bigoplus_{v>0}\,\bigl(N_v\oplus\overline{N}_v\bigr),
\end{split}
\end{equation}
where $N_v$ are complex vector bundles on which $S^1$ acts by sending $g\in S^1$ to $g^v$.
We also assume that (cf. \cite[(2.6)]{MR1870666})
\begin{equation}\label{s9}
V|_{F}=V_0^{\mathbb{R}}\oplus\bigoplus_{v>0}V_v\ ,\quad W|_{F}=\bigoplus_{v}W_v\ ,
\end{equation}
where $V_v$, $W_v$ are complex vector bundles on which $S^1$ acts by sending $g$ to $g^v$, and $V_0^{\mathbb{R}}$
is a real vector bundle on which $S^1$ acts as identity.

By \eqref{eg8}, as in \cite[(2.7)]{MR1870666} , there is a natural isomorphism
between the $\mathbb{Z}_2$-graded $C(TX)$-Clifford modules over $F$,
\begin{equation}\label{eg10}
S(TX,K_X)|_{F}\simeq S\Bigl(TY,K_X\otimes_{v>0}(\det N_v)^{-1}\Bigr)
\,\widehat{\otimes}\,\widehat{\bigotimes}_{v>0}\Lambda N_v\ .
\end{equation}

For a complex vector bundle $R$ over $F$, let $D^{Y}\otimes R$,
$D^{Y_\alpha}\otimes R$ be the twisted Spin$^c$ Dirac operator on
$S(TY,K_X\otimes_{v>0}(\det N_v)^{-1})\otimes R$ over $F$, $F_\alpha$ respectively.

We introduce the following locally constant functions on $F$ (cf. \cite[(2.8)]{MR1870666}),
\begin{equation}
\begin{split}\label{s14}
e(N)&=\sum_{v>0}v^2\dim N_v\ ,\quad d'(N)=\sum_{v>0}v\dim N_v\ ,
\\
e(V)&=\sum_{v>0}v^2\dim V_v\ ,\quad d'(V)=\sum_{v>0}v\dim V_v\ ,
\\
e(W)&=\sum_{v}v^2\dim W_v\ ,\quad d'(W)=\sum_{v}v\dim W_v\ .
\end{split}
\end{equation}

As in \cite[(2.9)]{MR1870666}, we write
\begin{equation}\label{eg17}
\begin{split}
L(N)&=\otimes_{v>0}(\det N_v)^v\ ,\quad L(V)=\otimes_{v>0}(\det V_v)^v\ ,
\\
L(W)&=\otimes_{v\neq0}(\det W_v)^v\ ,\quad
L=L(N)^{-1}\otimes L(V)\otimes L(W)^3\ .
\end{split}
\end{equation}
By using \eqref{hypothesis} and computing as in \cite[(2.10)-(2.11)]{MR1870666}, one knows
\begin{equation}\label{uy50}
c_1(L)=0\ ,\quad e(V)+3\cdot e(W)-e(N)=2e\ ,
\end{equation}
which means $L$ is a trivial complex line bundle over
each component $F_\alpha$ of $F$, and $S^1$ acts on $L$ by sending $g$ to $g^{2e}$, and $G_y$ acts on $L$
by sending $y$ to $y^{3d'(W)}$. From \cite[Lemma 2.1]{MR1870666}, we know that
$d'(W)\,{\rm mod}\,(N)$ is constant on each connected component of $M$.
Thus we can extend $L$
to a trivial complex line bundle over $M$, and we extend the $S^1$ action on it by
sending $g\in S^1$ on the canonical section $1$ of $L$ to $g^{2e}\cdot 1$, and $G_y$ acts
on $L$ by sending $y$ to $y^{3d'(W)}$.

In what follows, if $R(q)=\sum_{m\in\frac{1}{2}\mathbb{Z}}q^mR_m\in K_{S^1}(M)[[q^{1/2}]]$, we
will also denote $\ind(D^X\otimes R_m,h)$  by $\ind(D^X\otimes R(q),m,h)$. For $i=1,2,3,4$, set
\begin{equation}\label{s2}
R_{i1}=(K_W\otimes K_X^{-1})^{1/2}\otimes R_i(V)\otimes Q_1(W)\ .
\end{equation}

As in \cite[Proposition 2.1]{MR1870666}, by using Theorem \ref{local}
we first express the global equivariant family index
via the family indices on the fixed point set.

\begin{prop}\label{props3}
For $m\in \frac{1}{2}\mathbb{Z}$, $h\in \mathbb{Z}$, $1 \leq i\leq 4$, we have the
following identity in $K_{G_y}(B)$,
\begin{equation}
\begin{split}
&\ind\,\Bigl(D^X\otimes\otimes_{n=1}^{\infty}{\rm Sym}_{q^n}(TX)\otimes R_{i1},m,h\Bigr)
\\
&=\sum_{\alpha}(-1)^{\sum_{v>0}\dim N_v}
\ind\,\Bigl(D^{Y_\alpha}\otimes \otimes_{n=1}^{\infty}{\rm Sym}_{q^n}(TX|_F)\otimes R_{i1}\Bigr.
\\
\Bigl.&\hspace{90pt}\otimes{\rm Sym}\,(\oplus_{v>0}N_v)\otimes_{v>0}\det N_v,m,h\Bigr)\ .
\end{split}
\end{equation}
\end{prop}

To simplify the notations, we use the same convention as in \cite[pp. 945]{MR1870666}. For $n_0\in \mathbb{N}^*$, we define a number operator
$P$ on $K_{S^1}(M)[[q^{\frac{1}{n_0}}]]$ in the following way: if $R(q)=\oplus_{n\in\frac{1}{n_0}\mathbb{Z}}
R_nq^n\in K_{S^1}(M)[[q^{\frac{1}{n_0}}]]$, then $P$ acts on $R(q)$ by multiplication by $n$ on $R_n$. From
now on, we simply denote ${\rm Sym}_{q^n}(TX)$, $\Lambda_{q^n}(V)$ and $\Lambda_{q^n}(W)$
by ${\rm Sym}(TX_n)$, $\Lambda(V_n)$ and $\Lambda(W_n)$, respectively. In this way, $P$ acts on $TX_n$, $V_n$ and $W_n$
by multiplication by $n$, and the action of $P$ on ${\rm Sym}(TX_n)$, $\Lambda(V_{n})$ and $\Lambda(W_n)$ is naturally
induced by the corresponding action of $P$ on $TX_n$, $V_n$ and $W_n$. So the eigenspace of $P=n$ is just
given by the coefficient of $q^n$ of the corresponding element $R(q)$. For $R(q)=\oplus_{n\in\frac{1}{n_0}\mathbb{Z}}
R_nq^n\in K_{S^1}(M)[[q^{\frac{1}{n_0}}]]$, we will also denote
$\ind\bigl(D^{X}\otimes R_m,h\bigr)$ by $\ind\bigl(D^{X}\otimes R(q),m,h\bigr)$.

For $p\in \mathbb{N}$, we introduce the following elements in
$K_{S^1}(F)[[q]]$ (cf. \cite[(3.6)]{MR1870666}),
\begin{equation}\label{eg21}
\begin{split}
&\mathcal{F}_p(X)=\bigotimes_{n=1}^{\infty}{\rm Sym}(TY_n)\otimes\bigotimes_{v>0}\Bigl(\,\bigotimes_{n=1}^{\infty}
{\rm Sym}\,(N_{v,n})\bigotimes_{n>pv}{\rm Sym}\,(\overline{N}_{v,n})\,\Bigr)\ ,
\\
&\mathcal{F}_p'(X)=\bigotimes_{v>0}\bigotimes_{0\leq n\leq pv}\Bigl({\rm Sym}\,(N_{v,-n})\otimes\det N_v\Bigr)\ ,
\\
&\mathcal{F}^{-p}(X)=\mathcal{F}_p(X)\otimes\mathcal{F}_p'(X)\ .
\end{split}
\end{equation}
Then from \eqref{eg8}, over $F$, we have
\begin{equation}
\mathcal{F}^0(X)=\bigotimes_{n=1}^{\infty}{\rm Sym}_{q^n}(TX|_F)\otimes{\rm Sym}(\oplus_{v>0}N_v)\otimes_{v>0}\det N_v\ .
\end{equation}

We now state two intermediate results on the relations between the family indices on the fixed point set.
These two recursive formulas will be used in the next subsection to prove Theorem \ref{theorem}.

\begin{thm}\label{eg31}
{\rm (Compare with \cite[Theorem 2.3]{MR1870666})} For $1\leq i\leq 4$, $h$, $p\in \mathbb{Z}$, $p>0$,
$m\in \tfrac{1}{2}\mathbb{Z}$,  the following identity holds in $K_{G_y}(B)$,
\begin{equation}\label{eg32}
\begin{split}
&\sum_{\alpha}(-1)^{\sum_{v>0}\dim N_v}\ind\bigl(D^{Y_\alpha}\otimes\mathcal{F}^0(X)\otimes R_{i1},m,h\bigr)
\\
&=\sum_{\alpha}(-1)^{pd'(N)+\sum_{v>0}\dim N_v}\ind\bigl(D^{Y_\alpha}\otimes\mathcal{F}^{-p}(X)
\otimes R_{i1},\bigr.
\\
&\bigl.\hspace{75pt}
m+\frac{1}{2}p^2e(N)+\frac{p}{2}d'(N),h\bigr)\ .
\end{split}
\end{equation}
\end{thm}

\

The proof of Theorem \ref{eg31} will be given in Sections \ref{sec3.2}-\ref{sec3.4}.

\begin{thm}\label{eg29}
{\rm (Compare with \cite[Theorem 2.4]{MR1870666})} For each $\alpha$, $1\leq i\leq 4$,
$h$, $p\in \mathbb{Z}$, $p>0$, $m\in \tfrac{1}{2}\mathbb{Z}$,  the following identity holds in $K_{G_y}(B)$,
\begin{equation}
\begin{split}\label{eg30}
&\ind\Bigl(D^{Y_\alpha}\otimes {\cal F}^{-p}(X)\otimes R_{i1},m+\frac{1}{2}p^2e(N)+\frac{p}{2}d'(N),h\Bigr)
\\
&=(-1)^{pd'(W)}\ind\Bigl(D^{Y_\alpha}\otimes\mathcal{F}^{0}(X)\otimes R_{i1}\otimes L^{-p},m+ph+p^2e,h\Bigr)\ .
\end{split}
\end{equation}
\end{thm}

\

The proof of Theorem \ref{eg29} will be given in Section \ref{sec3.1}.

\subsection{A proof of Theorem \ref{theorem}}\label{sec2.4}

As $\tfrac{1}{2}p_1(3W-TX)_{S^1}\in H^*_{S^1}(X,\mathbb{Z})$ is well defined, one has the same identity as in \cite[(2.27)]{MR1870666},
\begin{equation}\label{s6}
d'(N)+d'(W)=0 \mod (2).
\end{equation}

From Proposition \ref{props3}, Theorems \ref{eg31}, \ref{eg29} and (\ref{s6}), for $1\leq i\leq 4$, $h, p\in\mathbb{Z}$, $p>0$, $m\in\tfrac{1}{2}\mathbb{Z}$, we get the following identity
(compare with \cite[(2.28)]{MR1870666}),
\begin{equation}
\begin{split}\label{s7}
&\ind\Bigl(D^{X}\otimes\bigotimes_{n=1}^{\infty}{\rm Sym}_{q^n}(TX)\otimes R_{i1},m,h\Bigr)
\\
&=\ind\Bigl(D^{X}\otimes\bigotimes_{n=1}^{\infty}{\rm Sym}_{q^n}(TX)\otimes R_{i1}\otimes L^{-p},m',h\Bigr).
\end{split}
\end{equation}
with
\begin{equation}\label{eg35}
m'=m+ph+p^2e.
\end{equation}

By (\ref{s3}) and (\ref{s2}), if $m<0$ or $m'<0$, then either side of (\ref{s7}) is identically zero,
which completes the proof of  Theorem \ref{theorem}. In fact,
\begin{enumerate}[{\rm (i)}]
\item Assume that $e=0$. Let $h\in \mathbb{Z}$, $m_0\in\tfrac{1}{2}\mathbb{Z}$,
$h\neq 0$ be fixed. If $h>0$, we take $m'=m_0$, then for $p$ large enough,
we get $m<0$ in (\ref{s7}). If $h<0$, we take $m=m_0$, then for $p$ large enough,
we get $m'<0$ in (\ref{s7}).

\item Assume that $e<0$. For $h\in\mathbb{Z}$, $m_0\in\tfrac{1}{2}\mathbb{Z}$,
we take $m=m_0$, then for $p$ large enough,
we get $m'<0$ in (\ref{s7}).
\end{enumerate}

The proof of  Theorem \ref{theorem} is completed.

\begin{rem}
We point out here that there is a $\mathbb{Z}/k$ version of Theorem \ref{theorem}, which is an analogue of
\cite[Theorem 4.4]{liuyu}. In fact,
by using the mod $k$ localization formula for
${\mathbb{Z}}/k$ circle actions on
${\mathbb{Z}}/k$ Spin$^c$ manifolds established in \cite[Theorem 2.7]{liuyu} (see also \cite[Theorem 2.1]{MR1993996} for the spin case),
our proof of Theorem \ref{theorem} can be applied to the case of $\mathbb{Z}/k$ manifolds with little modification.
\end{rem}

\begin{rem} {\rm (Compare with \cite[Remark 2.5]{MR1870666})} If M is connected, by \eqref{s7}, for $1\leq i\leq 4$, in $K_{G_y}(B)$,
we get
\begin{equation}
\begin{split}
&\ind\Bigl(D^{X}\otimes\bigotimes_{n=1}^{\infty}{\rm Sym}_{q^n}(TX)\otimes R_{i1}\Bigr)
\\
&\hspace{30pt}=\ind\Bigl(D^{X}\otimes\bigotimes_{n=1}^{\infty}{\rm Sym}_{q^n}(TX)\otimes R_{i1}\Bigr)\otimes [3d'(W)],
\end{split}
\end{equation}
where by $[3d'(W)]$ we mean the one dimensional complex vector space on which $y\in G_y$ acts by multiplication by $y^{3d'(W)}$.
In particular, if $B$ is a point, and $3d'(W)\neq 0 \mod(N)$, we get the vanishing theorem for String$^c$ manifolds analogue to
the result of \cite[$\S$\,10]{MR981372}.
\end{rem}

\section{Proofs of Theorems \ref{eg31} and \ref{eg29}}\label{sec3}

In this section, we prove those two intermediate results which are stated in Section \ref{sec2.3} and used
in Section \ref{sec2.4} to prove our main results.

This section is organized as follows. In Section \ref{sec3.1}, following \cite[Section 3.2]{MR1870666}, we prove Theorem \ref{eg29}.
In Section \ref{sec3.2}, we introduce the same refined shift operators as in \cite[Section 4.2]{MR1870666}.
In Section \ref{sec3.3}, we construct the twisted Spin$^c$ Dirac operator on $M(n_j)$, the fixed point
set of the naturally induced $\mathbb{Z}_{n_j}$-action on $M$. In Section \ref{sec3.4}, by applying the
$S^1$-equivariant index theorem in Section \ref{sec2.1}, we finally prove Theorem \ref{eg31}.

\subsection{A proof of Theorem \ref{eg29}}\label{sec3.1}

Let $H$ be the canonical basis of ${\rm Lie}(S^1)=\mathbb{R}$, i.e., $\exp(tH)=\exp(2\sqrt{-1}\pi t)$,
for $t\in \mathbb{R}$. On the fixed point $F$, let $\emph{\textbf{J}}_H$ denote the operator which computes
the weight of the $S^1$ action on $\Gamma(F,E|_{F})$ for any $S^1$-equivariant vector bundle $E$ over $M$.
Then $\emph{\textbf{J}}_H$ can be explicitly given by (cf. \cite[(3.2)]{MR2016198})
\begin{equation}
\emph{\textbf{J}}_H=\frac{1}{2\pi\sqrt{-1}}\mathscr{L}_H\big|_{\Gamma(F,E|_{F})}\ ,
\end{equation}
where $\mathscr{L}_H$ denotes the infinitesimal action of $H$ acting on $\Gamma(M,E)$.

Recall that the $\mathbb{Z}_2$-grading on $S(TX,K_X)\otimes_{n=1}^{\infty}{\rm Sym}(TX_n)$
(resp. $S(TY,K_X\otimes\otimes_{v>0}(\det N_v)^{-1})\otimes \mathcal{F}^{-p}(X)$)
is induced by the $\mathbb{Z}_2$-grading on $S(TX,K_X)$ (resp. $S(TY,K_X\otimes_{v>0}(\det N_v)^{-1})$). Write
\begin{align}\label{s10}
Q^1_{W}&=\bigotimes_{n=0}^\infty \Lambda(\overline{W}_n)\otimes \bigotimes_{n=1}^\infty\Lambda(W_n)\ ,\quad
Q^2_{W}=\bigotimes_{n\in \mathbb{N}+{1\over 2}}\Lambda(\overline{W}_n)
\otimes \bigotimes_{n\in \mathbb{N}+{1\over 2}}\Lambda(W_n)\ ,
\notag\\
F_V^1&=S(V)\otimes\bigotimes_{n=1}^{\infty}\Lambda(V_n)\ ,\quad
F_V^2=\bigotimes_{n\in\mathbb{N}+\frac{1}{2}}\Lambda(V_n)\ .
\end{align}
There are two natural $\mathbb{Z}_2$-gradings on $F_V^1$, $F_V^2$ (resp. $Q^1_W$, $Q^2_W$).
The first grading is induced by the $\mathbb{Z}_2$-grading
of $S(V)$ and the forms of homogeneous degrees in $\otimes_{n=1}^{\infty}\Lambda(V_n)$,
$\otimes_{n\in\mathbb{N}+\frac{1}{2}}\Lambda(V_n)$ (resp. $Q^2_W$). We define $\tau_e|_{F_V^{i\pm}}=\pm 1$ ($i=1,2$)
(resp. $\tau_e|_{Q^{2\pm}_W}=\pm1$) to be the involution
defined by this $\mathbb{Z}_2$-grading. The second grading is the one for which $F_V^i$ and $Q_W^i$ ($i=1,2$) are purely even,
i.e., $F_V^{i+}=F_V^i$, $Q_W^{i+}=Q_W^i$. We denote by $\tau_s={\rm id}$ the involution defined by this $\mathbb{Z}_2$-grading.
Set $Q(W)=Q^1_{W}\otimes Q^2_{W}\otimes Q^2_{W}$. We will denote by $\tau_1$ the $\mathbb{Z}_2$-grading
on $Q(W)$ defined by
\begin{equation}
(Q(W),\tau_1)=(Q^1_W,\tau_s)\widehat\otimes (Q^2_W,\tau_e)\widehat\otimes (Q^2_W,\tau_s).
\end{equation}
Then the coefficient of $q^n$ ($n\in\frac{1}{2}\mathbb{Z}$) in (\ref{s3}) of $R_1(V)$
(resp.  $R_2(V)$, $R_3(V)$, $R_{4}(V)$, $Q_1(W)$) is exactly the $\mathbb{Z}_2$-graded vector subbundle of $(F_V^1,\tau_s)$
(resp. $(F_V^1,\tau_e)$, $(F_V^2,\tau_e)$, $(F_V^2,\tau_s)$, $(Q(W),\tau_1)$), on which $P$ acts by multiplication by $n$.

Furthermore, we denote by $\tau_e$ (resp. $\tau_s$) the $\mathbb{Z}_2$-grading on $S(TX,K_X)
\otimes\otimes_{n=1}^{\infty}{\rm Sym}(TX_n)\otimes F_V^i$ ($i=1$, $2$) induced by
the above $\mathbb{Z}_2$-gradings. We will denote by $\tau_{e1}$ (resp. $\tau_{s1}$)
the $\mathbb{Z}_2$-grading on $S(TX,K_X)\otimes\otimes_{n=1}^{\infty}{\rm Sym}(TX_n)\otimes F_V^i\otimes Q(W)$ ($i=1$, $2$) defined by
\begin{equation}\label{eg13}
\tau_{e 1}=\tau_e\widehat{\otimes}\tau_1\ ,\quad \tau_{s 1}=\tau_s\widehat{\otimes}\tau_1\ .
\end{equation}
We still denote by $\tau_{e1}$ (resp. $\tau_{s1}$) the $\mathbb{Z}_2$-grading on $S(TY,K_X\otimes_{v>0}(\det N_v)^{-1})\otimes \mathcal{F}^{-p}(X)\otimes F_V^i\otimes Q(W)$ ($i=1$, $2$) which is induced as in \eqref{eg13}.

By (\ref{s9}), as in (\ref{eg10}), there is a natural isomorphism between the
$\mathbb{Z}_2$-graded $C(V)$-Clifford modules over $F$,
\begin{equation}\label{eg14}
S(V)|_{F}\simeq S\Bigl(V_0^\mathbb{R},\otimes_{v>0}(\det V_v)^{-1}\Bigr)\otimes\widehat{\bigotimes}_{v>0}\Lambda V_v\ .
\end{equation}

Let $V_0=V_0^{\mathbb{R}}\otimes_{\mathbb{R}}\mathbb{C}$. Using \eqref{s9} and \eqref{eg14},
we rewrite \eqref{s10} on the fixed point set $F$ as follows,
\begin{equation}\label{eg15}
\begin{split}
Q^1_W&=\bigotimes_{n=0}^{\infty}\Lambda(\oplus_{v}\overline{W}_{v,n})
\otimes\bigotimes_{n=1}^{\infty}\Lambda(\oplus_{v}W_{v,n})\ ,
\\
Q^2_W&=\bigotimes_{n\in \mathbb{N}+{1\over 2}}\Lambda(\oplus_{v}\overline{W}_{v,n})
\otimes\bigotimes_{n\in \mathbb{N}+{1\over 2}}\Lambda(\oplus_{v}W_{v,n})\ ,
\\
F_V^1&=\bigotimes_{n=1}^{\infty}\Lambda\Bigl(V_{0,n}\oplus\oplus_{v>0}(V_{v,n}\oplus\overline{V}_{v,n})\Bigr)
\\
&\hspace{70pt}\otimes S\Bigl(V_0^\mathbb{R},\otimes_{v>0}(\det V_v)^{-1}\Bigr)\otimes_{v>0}\Lambda V_{v,0}\ ,
\\
F_V^2&=\bigotimes_{n\in\mathbb{N}+\frac{1}{2}}\Lambda\Bigl(V_{0,n}\oplus\oplus_{v>0}(V_{v,n}\oplus\overline{V}_{v,n})\Bigr)\ .
\end{split}
\end{equation}

We can reformulate Theorem \ref{eg29} as follows.
\begin{thm}\label{s15}
For each $\alpha$, $h$, $p\in \mathbb{Z}$, $p>0$,
$m\in \tfrac{1}{2}\mathbb{Z}$, for $i=1,2$, $\tau=\tau_{e1}$ or $\tau_{s1}$,
the following identity holds in $K_{G_y}(B)$,
\begin{equation}
\begin{split}
&\ind_\tau\Big(D^{Y_\alpha}\otimes(K_W\otimes K_X^{-1})^{1/2}\otimes \mathcal{F}^{-p}(X)\otimes F^i_V\otimes Q(W),
\\
&\hspace{100pt}m+{1\over 2}p^2e(N)+{1\over 2}pd'(N),h\Big)
\\
&=(-1)^{pd'(W)}\ind_\tau\Big(D^{Y_\alpha}\otimes(K_W\otimes K_X^{-1})^{1/2}\otimes \mathcal{F}^0(X)\otimes F^i_V
\\
&\hspace{105pt}\otimes Q(W)\otimes L^{-p},m+ph+p^2e,h\Big)\ .
\end{split}
\end{equation}
\end{thm}

Following \cite{MR998662} in spirit, we introduce the same shift operators as in \cite[(3.9)]{MR1870666}.
For $p\in \mathbb{N}$, we set
\begin{equation}\label{eg16}
\begin{split}
&r_*:N_{v,n}\rightarrow N_{v,n+pv}\ ,\quad r_*:\overline{N}_{v,n}\rightarrow \overline{N}_{v,n-pv}\ ,
\\
&r_*:V_{v,n}\rightarrow V_{v,n+pv}\ ,\quad \ r_*:\overline{V}_{v,n}\rightarrow \overline{V}_{v,n-pv}\ ,
\\
&r_*:W_{v,n}\rightarrow W_{v,n+pv}\ ,\quad r_*:\overline{W}_{v,n}\rightarrow \overline{W}_{v,n-pv}\ .
\end{split}
\end{equation}

\begin{prop}\label{s8}
For $p\in \mathbb{Z}$, $p>0$, $i=1,2$, there are natural isomorphisms of vector bundles over $F$,
\begin{equation}\label{eg23}
r_*(\mathcal{F}^{-p}(X))\simeq\mathcal{F}^{0}(X)\otimes L(N)^{p}\ ,\quad
r_*(F_V^i)\simeq F_V^i\otimes L(V)^{-p}\ .
\end{equation}
For any $p\in \mathbb{Z}$, $p>0$, $i=1,2$,
there are natural $G_y\times S^1$-equivariant isomorphisms of vector bundles over $F$,
\begin{equation}\label{eg24}
r_*(Q^i_W)\simeq Q^i_W\otimes L(W)^{-p}\ .
\end{equation}
In particular, one gets the $G_y\times S^1$-equivariant bundle isomorphism
\begin{equation}\label{s12}
r_*(Q(W))\simeq Q(W)\otimes L(W)^{-3p}\ .
\end{equation}
\end{prop}
{\bf Proof}\hspace{4mm}By \cite[Proposition 3.1]{MR1870666}, only the $i=2$ case in \eqref{eg24} needs to be proved.

Using \cite[(3.14)-(3.16)]{MR1870666}, we have the natural $G_y\times S^1$-equivariant isomorphisms of vector bundles over $F$,
\begin{equation}
\begin{split}\label{s11}
\bigotimes_{n\in\mathbb{N}+{1\over 2},v>0,\atop 0<n<pv}\Lambda^{i_n}(\overline{W}_{v,n-pv})
&\simeq\bigotimes_{n\in\mathbb{N}+{1\over 2},v>0,\atop 0<n<pv}\Lambda^{\dim W_v-i_n}({W}_{v,-n+pv})
\\
&\hspace{90pt}\otimes\bigotimes_{v>0}(\det \overline{W}_v)^{pv}\ ,
\\
\bigotimes_{n\in\mathbb{N}+{1\over 2},v<0,\atop 0<n<-pv}\Lambda^{i'_n}(W_{v,n+pv})
&\simeq\bigotimes_{n\in\mathbb{N}+{1\over 2},v<0,\atop 0<n<-pv}\Lambda^{\dim W_v-i'_n}(\overline{W}_{v,-n-pv})
\\
&\hspace{90pt}\otimes\bigotimes_{v>0}(\det {W}_v)^{-pv}\ .
\end{split}
\end{equation}
From \eqref{eg17} and \eqref{s11}, we get \eqref{eg24} for the $i=2$ case.

The proof of Proposition \ref{s8} is completed.

\

The following proposition, which is an analogue of \cite[Proposition 3.2]{MR1870666}, is deduced from Proposition \ref{s8}.
\begin{prop}\label{eg25}
For $p\in \mathbb{Z}$, $p>0$, $i=1$, $2$, the $G_y$-equivariant isomorphism of vector bundles over $F$
induced by {\rm (\ref{eg23})}, {\rm (\ref{s12})},
\begin{equation}\label{eg26}
\begin{split}
r_*\,:&\ S(TY,K_X\otimes_{v>0}(\det N_v)^{-1})\otimes(K_W\otimes K^{-1}_X)^{1/2}
\\
&\hspace{50pt}\otimes\mathcal{F}^{-p}(X)
\otimes F^i_V\otimes Q(W)
\\
&\longrightarrow\ S(TY,K_X\otimes_{v>0}(\det N_v)^{-1})\otimes(K_W\otimes K^{-1}_X)^{1/2}
\\
&\hspace{50pt}\otimes\mathcal{F}^{0}(X)
\otimes F^i_V\otimes Q(W)\otimes L^{-p}\ ,
\end{split}
\end{equation}
verifies the following identities,
\begin{equation}\label{eg27}
\begin{split}
r_*^{-1}\cdot \textbf{J}_H\cdot r_*&=\textbf{J}_H\ ,
\\
r_*^{-1}\cdot P\cdot r_*&=P+p\textbf{J}_H+p^2e-\frac{1}{2}p^2e(N)-\frac{p}{2}d'(N)\ .
\end{split}
\end{equation}
For the $\mathbb{Z}_2$-gradings, we have
\begin{equation}\label{eg28}
r_*^{-1}\tau_{e}r_*=\tau_{e}\ ,
\quad r_*^{-1}\tau_{s}r_*=\tau_{s}\ ,
\quad r_*^{-1}\tau_{1}r_*=(-1)^{pd'(W)}\tau_{1}\ .
\end{equation}
\end{prop}
{\bf Proof}\hspace{4mm}By the proof of  \cite[Proposition 3.2]{MR1870666}, we need to compute the action of
$r_*^{-1}\cdot P\cdot r_*$ on $\bigotimes_{n\in\mathbb{N}+{1\over 2},v>0,\atop 0<n<pv}
\Lambda^{i_n}(\overline{W}_{v,n})
\bigotimes_{n\in\mathbb{N}+{1\over 2},v<0,\atop 0<n<-pv}
\Lambda^{i'_n}(W_{v,n})$. In fact, by \eqref{s11},
\begin{equation}
\begin{split}\label{s13}
r_*^{-1}\cdot P\cdot r_*&=\sum_{n\in\mathbb{N}+{1\over 2},v>0,\atop 0<n<pv}(\dim W_v-i_n)(-n+pv)
\\
&\hspace{60pt}+\sum_{n\in\mathbb{N}+{1\over 2},v<0,\atop 0<n<-pv}(\dim W_v-i'_n)(-n-pv)
\\
&=P+p\emph{\textbf{J}}_H+{1\over 2}p^2e(W)\ .
\end{split}
\end{equation}
By \cite[(3.21)-(3.23)]{MR1870666}, \eqref{s14}-\eqref{uy50} and \eqref{s13}, we deduce the second line of \eqref{eg27}.
The first line of \eqref{eg27} is obvious.

Consider the $\mathbb{Z}_2$-gradings. The first two identities of \eqref{eg28} were proved in \cite[(3.18)]{MR2016198}.
$\tau_1$ changes only on
$\bigotimes_{n\in\mathbb{N}+{1\over 2},v>0,\atop 0<n<pv}
\Lambda^{i_n}(\overline{W}_{v,n})
\bigotimes_{n\in\mathbb{N}+{1\over 2},v<0,\atop 0<n<-pv}
\Lambda^{i'_n}(W_{v,n})$. From \eqref{s14} and \eqref{s11}, we get the third identity of \eqref{eg28}.

The proof of Proposition \ref{eg25} is completed.

\

Theorem \ref{s15} is a direct consequence of Proposition \ref{eg25}. The proof of Theorem \ref{eg29} is completed.

The rest of this section is devoted to a proof of Theorem \ref{eg31}.

\subsection{The refined shift operators}\label{sec3.2}

We first introduce a partition of $[0,1]$ as in \cite[pp. 942--943]{MR1870666}. Set
$J=\bigl\{v\in \mathbb{N}\,\big| \text{ there exists } \alpha \text{ such that } N_v\neq0 \text{ on } F_\alpha\bigr\}$ and
\begin{equation}\label{eg39}
\Phi=\bigl\{\beta\in (0,1]\,\big| \text{ there exists } v\in J \text{ such that } \beta v\in\mathbb{Z}\bigr\}\ .
\end{equation}
We order the elements in $\Phi$ so that $\Phi=\bigl\{\beta_i\,\big|\,1\leq i\leq J_0, J_0\in\mathbb{N}\text{ and }
\beta_i<\beta_{i+1}\bigr\}$. Then for any integer $1\leq i\leq J_0$, there exist $p_i$, $n_i\in\mathbb{N}$,
$0<p_i\leq n_i$, with $(p_i,n_i)=1$ such that
\begin{equation}\label{eg40}
\beta_i=p_i/n_i\ .
\end{equation}
Clearly, $\beta_{J_0}=1$. We also set $p_0=0$ and $\beta_0=0$.

For $0\leq j\leq J_0$, $p\in\mathbb{N}^*$, we write
\begin{equation}\label{eg41}
\begin{split}
I_j^p&=\Bigl\{(v,n)\in\mathbb{N}\times\mathbb{N}\,\Big|\,v\in J,\ (p-1)v<n\leq pv,
\ \frac{n}{v}=p-1+\frac{p_j}{n_j}\Bigr\}\ ,
\\
\overline{I}_j^p&=\Bigl\{(v,n)\in\mathbb{N}\times\mathbb{N}\,\Big|\,v\in J,\ (p-1)v<n\leq pv,\
\frac{n}{v}>p-1+\frac{p_j}{n_j}\Bigr\}\ .
\end{split}
\end{equation}
Clearly, $I_0^p=\emptyset$, the empty set. We define $\mathcal{F}_{p,j}(X)$ as in \cite[(2.21)]{MR1870666},
which are analogous with (\ref{eg21}). More specifically, we set
\begin{align}
\mathcal{F}_{p,j}(X)&=\bigotimes_{n=1}^{\infty}{\rm Sym}\,(TY_n)\otimes\bigotimes_{v>0}\Bigl(\,\bigotimes_{n=1}^{\infty}
{\rm Sym}\,(N_{v,n})\otimes\bigotimes_{n> (p-1)v +\frac{p_j}{n_j}v}{\rm Sym}\,(\overline{N}_{v,n})\,\Bigr)
\notag\\
&\hspace{65pt}\otimes\bigotimes_{v>0,\atop 0\leq n\leq (p-1)v+\bigl[\frac{p_j}{n_j}v\bigr]}
\Bigl({\rm Sym}\,(N_{v,-n})\otimes\det N_v\Bigr)
\notag\\
&=\mathcal{F}_{p}(X)\otimes\mathcal{F}_{p-1}'(X)\otimes\bigotimes_{(v,n)\in\overline{I}_j^p}{\rm Sym}\,(\overline{N}_{v,n})
\\
&\hspace{65pt}\otimes\bigotimes_{(v,n)\in\cup_{i=0}^jI_i^p}
\Bigl({\rm Sym}\,(N_{v,-n})\otimes\det N_v\Bigr)\ ,\notag
\end{align}
where we use the notation that for $s\in \mathbb{R}$, $[s]$ denotes the greatest integer which is less than or equal to $s$. Then
\begin{equation}
\mathcal{F}_{p,0}(X)=\mathcal{F}^{-p+1}(X)\ ,\quad
\mathcal{F}_{p,J_0}(X)=\mathcal{F}^{-p}(X)\ .
\end{equation}

From the construction of $\beta_i$, we know that for $v\in J$, there is no integer in
$\bigl(\tfrac{p_{j-1}}{n_{j-1}}v,\tfrac{p_j}{n_j}v\bigr)$. Furthermore (cf. \cite[(4.24)]{MR1870666}),
\begin{equation}
\begin{split}\label{pj}
\Bigl[\frac{p_{j-1}}{n_{j-1}}v\Bigr]&=\Bigl[\frac{p_{j}}{n_{j}}v\Bigr]-1\ \text{  if  }\ v\equiv0\mod(n_j)\ ,
\\
\Bigl[\frac{p_{j-1}}{n_{j-1}}v\Bigr]&=\Bigl[\frac{p_{j}}{n_{j}}v\Bigr]\ \text{ if }\ v\not \equiv0\mod(n_j)\ .
\end{split}
\end{equation}

We use the same shift operators $r_{j*}$, $1\leq j\leq J_0$ as in \cite[(4.21)]{MR1870666}, which
refine the shift operator $r_{*}$ defined in (\ref{eg16}). For $p\in\mathbb{N}\backslash\{0\}$, set
\begin{equation}\label{eg43}
\begin{split}
&r_{j*}:N_{v,n}\rightarrow N_{v,n+(p-1)v+p_jv/n_j}\ ,\quad
r_{j*}:\overline{N}_{v,n}\rightarrow \overline{N}_{v,n-(p-1)v-p_jv/n_j}\ ,
\\
&r_{j*}:V_{v,n}\rightarrow V_{v,n+(p-1)v+p_jv/n_j}\ ,
\quad \ r_{j*}:\overline{V}_{v,n}\rightarrow \overline{V}_{v,n-(p-1)v-p_jv/n_j}\ ,
\\
&r_{j*}:W_{v,n}\rightarrow W_{v,n+(p-1)v+p_jv/n_j}\ ,
\quad r_{j*}:\overline{W}_{v,n}\rightarrow \overline{W}_{v,n-(p-1)v-p_jv/n_j}\ .
\end{split}
\end{equation}

For $1\leq j\leq J_0$,  we define $\mathcal{F}(\beta_j)$, $F_V^1(\beta_j)$, $F_V^2(\beta_j)$,
$Q^1_W(\beta_j)$ and $Q^2_W(\beta_j)$ over $F$ as follows (compare with \cite[(4.13)]{MR1870666}).
\begin{equation*}
\begin{split}
&\mathcal{F}(\beta_j)=\bigotimes_{0<n\in\mathbb{Z}}{\rm Sym}\,(TY_n)\otimes\bigotimes_{v>0,\atop v\equiv0,\frac{n_j}{2}
\,{\rm mod}(n_j) }\bigotimes_{0<n\in\mathbb{Z}+\frac{p_j}{n_j}v}{\rm Sym}\,(N_{v,n}\oplus\overline{N}_{v,n})
\\
&\quad\otimes\bigotimes_{0<v'<n_j/2}{\rm Sym}\,\Biggl(\ \bigoplus_{v\equiv v',-v'\ {\rm mod}(n_j)}\Bigl(\bigoplus_
{0<n\in\mathbb{Z}+\frac{p_j}{n_j}v}N_{v,n}\bigoplus_{0<n\in\mathbb{Z}-\frac{p_j}{n_j}v}
\overline{N}_{v,n}\Bigr)\ \Biggr)\ ,
\end{split}
\end{equation*}
\begin{equation*}
\begin{split}
F_V^1(\beta_j)&=\Lambda\Bigg(\ \bigoplus_{0<n\in\mathbb{Z}}V_{0,n}\bigoplus_{v>0,\atop v\equiv0,\frac{n_j}{2}
\,{\rm mod}(n_j)}\Bigl(\bigoplus_
{0<n\in\mathbb{Z}+\frac{p_j}{n_j}v}V_{v,n}\bigoplus_{0<n\in\mathbb{Z}-\frac{p_j}{n_j}v}
\overline{V}_{v,n}\Bigr)
\\
&\bigoplus_{0<v'<n_j/2}\Bigl(\ \bigoplus_{v\equiv v',-v'\,{\rm mod}(n_j)}\Bigl(\bigoplus_
{0<n\in\mathbb{Z}+\frac{p_j}{n_j}v}V_{v,n}\bigoplus_{0<n\in\mathbb{Z}-\frac{p_j}{n_j}v}
\overline{V}_{v,n}\Bigr)\Bigr)\Bigg)\ ,
\end{split}
\end{equation*}
\begin{equation*}
\begin{split}
F_V^2(\beta_j)\!&=\!\Lambda\Bigg(\ \bigoplus_{0<n\in\mathbb{Z}+\frac{1}{2}}V_{0,n}
\bigoplus_{v>0,\atop v\equiv0,\frac{n_j}{2}\,{\rm mod}(n_j)}\Bigl(\bigoplus_{\negthickspace 0<n\in\mathbb{Z}+\frac{p_j}{n_j}v+\frac{1}{2}}V_{v,n}
\negthickspace\bigoplus_{\negthickspace 0<n\in\mathbb{Z}-\frac{p_j}{n_j}v+\frac{1}{2}}
\overline{V}_{v,n}\!\Bigr)
\\
&\bigoplus_{0<v'<n_j/2}\Bigl(\ \bigoplus_{v\equiv v',-v'\,{\rm mod}(n_j)}\Bigl(\bigoplus_
{0<n\in\mathbb{Z}+\frac{p_j}{n_j}v+\frac{1}{2}}V_{v,n}\bigoplus_{0<n\in\mathbb{Z}-\frac{p_j}{n_j}v+\frac{1}{2}}
\overline{V}_{v,n}\Bigr)\Bigr)\Bigg)\ ,
\end{split}
\end{equation*}
\begin{align}
\label{s18}
Q^1_W(\beta_j)&=\Lambda\Biggl(\ \bigoplus_v\ \Bigl(\bigoplus_
{0<n\in\mathbb{Z}+\frac{p_j}{n_j}v}W_{v,n}\bigoplus_{0\leq n\in\mathbb{Z}-\frac{p_j}{n_j}v}
\overline{W}_{v,n}\ \Bigr)\Biggr)\ ,
\notag\\
Q^2_W(\beta_j)&=\Lambda\Biggl(\ \bigoplus_v\ \Bigl(\bigoplus_
{0<n\in\mathbb{Z}+\frac{p_j}{n_j}v+{1\over 2}}W_{v,n}\bigoplus_{0<n\in\mathbb{Z}-\frac{p_j}{n_j}v+{1\over 2}}
\overline{W}_{v,n}\ \Bigr)\Biggr).
\end{align}

Using \eqref{pj}, \eqref{s18} and computing directly, we get an analogue of
\cite[Proposition 4.1]{MR1870666} which refines Proposition \ref{s8}.
\begin{prop}\label{eg44}
For $p\in \mathbb{Z}$, $p>0$, $1\leq j\leq J_0$, there are natural isomorphisms of vector bundles over $F$,
\begin{equation*}
\begin{split}
r_{j*}(\mathcal{F}_{p,j-1}(X))&\simeq\mathcal{F}(\beta_j)\otimes\bigotimes_{v>0,\ v\equiv0\ {\rm mod}\,(n_j)}
{\rm Sym}\,(\overline{N}_{v,0})
\\
&\otimes\bigotimes_{v>0}\ (\det N_v)^{\bigl[\frac{p_j}{n_j}v\bigr]+(p-1)v+1}
\otimes\bigotimes_{v>0,\ v\equiv0\ {\rm mod}\,(n_j)}(\det N_v)^{-1}\ ,
\\
r_{j*}(\mathcal{F}_{p,j}(X))&\simeq\ \mathcal{F}(\beta_j)\otimes\bigotimes_{v>0,\ v\equiv0\ {\rm mod}\,(n_j)}
{\rm Sym}(N_{v,0})
\\
&\hspace{110pt}
\otimes\bigotimes_{v>0}\ (\det N_v)^{\bigl[\frac{p_j}{n_j}v\bigr]+(p-1)v+1}\ ,
\end{split}
\end{equation*}
\begin{equation*}
\begin{split}
r_{j*}(F_V^1)& \simeq S\Bigl(V_0^{\mathbb{R}}\,,\,\otimes_{v>0}(\det V_v)^{-1}\Bigr)\otimes F_V^1(\beta_j)
\\
&\hspace{40pt}\otimes\bigotimes_{v>0,\atop v\equiv0\ {\rm mod}\,(n_j)}\Lambda(V_{v,0})\otimes
\bigotimes_{v>0}(\det \overline{V}_v)^{\bigl[\frac{p_j}{n_j}v\bigr]+(p-1)v}\ ,
\\
r_{j*}(F_V^2)& \simeq F_V^2(\beta_j)\otimes\bigotimes_{v>0,\atop v\equiv\frac{n_j}{2}\ {\rm mod}\,(n_j)}
\Lambda(V_{v,0})\otimes\bigotimes_{v>0}(\det \overline{V}_v)^{\bigl[\frac{p_j}{n_j}v+\frac{1}{2}\bigr]+(p-1)v}\ .
\end{split}
\end{equation*}
For $p\in \mathbb{Z}$, $p>0$, $1\leq j\leq J_0$, there are natural $G_y\times S^1$-equivariant isomorphisms of vector bundles over $F$,
\begin{equation}
\begin{split}\label{s16}
r_{j*}(Q_W^1)& \simeq Q_W^1(\beta_j)\otimes\bigotimes_{v>0,\atop v\equiv0\ {\rm mod}\,(n_j)}\det W_{v}\ \otimes
\\
&\quad
\bigotimes_{v>0}(\det \overline{W}_v)^{\bigl[\frac{p_j}{n_j}v\bigr]+(p-1)v+1}
\bigotimes_{v<0}(\det {W}_v)^{\bigl[-\frac{p_j}{n_j}v\bigr]-(p-1)v}\ ,
\\
r_{j*}(Q_W^2)& \simeq Q_W^2(\beta_j)\otimes\bigotimes_{v>0,\atop v\equiv \frac{n_j}{2}\ {\rm mod}\,(n_j)}\Lambda(W_{v,0})
\otimes\bigotimes_{v<0,\atop v\equiv \frac{n_j}{2}\ {\rm mod}\,(n_j)}\Lambda(\overline{W}_{v,0})
\\
&\bigotimes_{v>0}(\det \overline{W}_v)^{\bigl[\frac{p_j}{n_j}v+\frac{1}{2}\bigr]+(p-1)v}
\bigotimes_{v<0}(\det {W}_v)^{\bigl[-\frac{p_j}{n_j}v+\frac{1}{2}\bigr]-(p-1)v}\ .
\end{split}
\end{equation}
\end{prop}
{\bf Proof}\hspace{4mm}By \cite[Proposition 4.1]{MR1870666}, we need only to prove the second isomorphism in \eqref{s16}.
In fact, using  \cite[(3.14)]{MR1870666}, we have the natural $G_y\times S^1$-equivariant isomorphisms of vector bundles over $F$,
\begin{equation}
\begin{split}\label{s17}
&\bigotimes_{0<n\in\mathbb{Z}+{1\over 2},\,v>0,\atop n-(p-1)v-\frac{p_j}{n_j}v\leq 0}
\Lambda^{i_n}\Bigl(\overline{W}_{v,n-(p-1)v-\frac{p_j}{n_j}v}\Bigr)\simeq
\bigotimes_{v>0}(\det \overline{W}_v)^{\bigl[\frac{p_j}{n_j}v+\frac{1}{2}\bigr]+(p-1)v}
\\
&\hspace{80pt}\otimes\bigotimes_{0<n\in\mathbb{Z}+{1\over 2},\,v>0,\atop n-(p-1)v-\frac{p_j}{n_j}v\leq 0}
\Lambda^{\dim W_v-i_n}\Bigl({W}_{v,-n+(p-1)v+\frac{p_j}{n_j}v}\Bigr)\ ,
\\
&\bigotimes_{0<n\in\mathbb{Z}+{1\over 2},\,v<0,\atop n+(p-1)v+\frac{p_j}{n_j}v\leq 0}
\Lambda^{i'_n}\Bigl({W}_{v,n+(p-1)v+\frac{p_j}{n_j}v}\Bigr)\simeq\bigotimes_{v<0}
(\det {W}_v)^{\bigl[-\frac{p_j}{n_j}v+\frac{1}{2}\bigr]-(p-1)v}
\\
&\hspace{80pt}\otimes\bigotimes_{0<n\in\mathbb{Z}+{1\over 2},\,v<0,\atop n+(p-1)v+\frac{p_j}{n_j}v\leq 0}
\Lambda^{\dim W_v-i'_n}\Bigl(\overline{W}_{v,-n-(p-1)v-\frac{p_j}{n_j}v}\Bigr)\ ,
\end{split}
\end{equation}
From \eqref{s18} and \eqref{s17}, we get the second isomorphism in \eqref{s16}.

The proof of Proposition \ref{eg44} is completed.

\subsection{The Spin$^c$ Dirac operators on $M(n_j)$}\label{sec3.3}

Recall that there is a nontrivial circle action on $M$ which can be lifted to
the circle actions on $V$ and $W$.

For $n\in \mathbb{N}\backslash\{0\}$, let $\mathbb{Z}_n\subset S^1$
denote the cyclic subgroup of order $n$. Let $M(n_j)$ be the fixed point set of the induced $\mathbb{Z}_{n_j}$
action on $M$. Then $\pi:M(n_j)\rightarrow B$ is a fibration with compact fiber $X(n_j)$.
Let $N(n_j)\rightarrow M(n_j)$ be the normal bundle to $M(n_j)$ in $M$. As in \cite[pp. 151]{MR954493}
(see also \cite[Section 4.1]{MR1870666}, \cite[Section 4.1]{MR2016198} or \cite{MR998662}),
we see that $N(n_j)$ and $V$ can be decomposed, as real vector bundles over $M(n_j)$, into
\begin{equation}\label{eg48}
\begin{split}
N(n_j)&=\bigoplus_{0<v<n_j/2}N(n_j)_v\oplus N(n_j)_{n_j/2}^{\mathbb{R}}\ ,
\\
V|_{M(n_j)}&= V(n_j)_{0}^{\mathbb{R}}\oplus\bigoplus_{0<v<n_j/2}V(n_j)_v\oplus V(n_j)_{n_j/2}^{\mathbb{R}}\ ,
\end{split}
\end{equation}
where $V(n_j)_{0}^{\mathbb{R}}$ is the real vector bundle on which $\mathbb{Z}_{n_j}$ acts by identity,
and $N(n_j)_{n_j/2}^{\mathbb{R}}$ (resp. $V(n_j)_{n_j/2}^{\mathbb{R}}$) is defined to be zero if $n_j$ is odd.
Moreover, for $0<v<n_j/2$, $N(n_j)_v$ (resp. $V(n_j)_v$) admits a
unique complex structure such that $N(n_j)_v$ (resp. $V(n_j)_v$) becomes a
complex vector bundle on which $g\in \mathbb{Z}_{n_j}$ acts by $g^v$. We also denote by
$V(n_j)_0$, $V(n_j)_{n_j/2}$ and $N(n_j)_{n_j/2}$ the corresponding complexification of
$V(n_j)_{0}^{\mathbb{R}}$, $V(n_j)_{n_j/2}^{\mathbb{R}}$ and $N(n_j)_{n_j/2}^{\mathbb{R}}$.

Similarly, we also have the following $\mathbb{Z}_{n_j}$-equivariant decomposition of $W$ over $M(n_j)$ into complex vector bundles,
\begin{equation}\label{eg49}
W|_{M(n_j)}=\bigoplus_{0\leq v<n_j}W(n_j)_v\ ,
\end{equation}
where for $0\leq v<n_j$, $g\in \mathbb{Z}_{n_j}$ acts on $W(n_j)_v$ by sending $g$ to $g^v$.

By \cite[Lemma 4.1]{MR1870666} (which generalizes \cite[Lemmas 9.4 and 10.1]{MR954493} and \cite[Lemma 5.1]{MR998662}),
we know that the vector bundles $TX(n_j)$ and $V(n_j)_{0}^{\mathbb{R}}$ are orientable and even dimensional. Thus $N(n_j)$ is orientable
over $M(n_j)$. By (\ref{eg48}), $V(n_j)_{n_j/2}^{\mathbb{R}}$ and $N(n_j)_{n_j/2}^{\mathbb{R}}$ are also
orientable and even dimensional. In what follows, we fix the orientations of $N(n_j)_{n_j/2}^{\mathbb{R}}$
and $V(n_j)_{n_j/2}^{\mathbb{R}}$. Then $TX(n_j)$ and $V(n_j)_0^{\mathbb{R}}$ are naturally oriented by
(\ref{eg48}) and the orientations of $TX$, $V$, $N(n_j)_{n_j/2}^{\mathbb{R}}$ and
$V(n_j)_{n_j/2}^{\mathbb{R}}$. Let $W(n_j)^\mathbb{R}_{{n_j}/{2}}$ be the underlying real vector bundle of
$W(n_j)_{{n_j}/{2}}$, which are canonically oriented by its complex structure.

By (\ref{eg8}), (\ref{s9}), (\ref{eg48}) and (\ref{eg49}), we get
the following identifications of complex vector bundles over $F$ (cf. \cite[(4.9) and (4.12)]{MR1870666}),
for $0<v\leq n_j/2$,
\begin{equation}
\begin{split}\label{uy2}
N(n_j)_v\big|_{F}&=\bigoplus_{v'>0,\,v'\equiv v\,{\rm mod}(n_j)}N_{v'}\ \oplus
\bigoplus_{v'>0,\,v'\equiv -v\,{\rm mod}(n_j)}\overline{N}_{v'}\ \ ,
\\
V(n_j)_v\big|_{F}&=\bigoplus_{v'>0,\,v'\equiv v\,{\rm mod}(n_j)}V_{v'}\ \oplus
\bigoplus_{v'>0,\,v'\equiv -v\,{\rm mod}(n_j)}\overline{V}_{v'}\ \ ,
\end{split}
\end{equation}
for $0\leq v<n_j$,
\begin{equation}\label{uy3}
W(n_j)_v\big|_{F}=\bigoplus_{v'\equiv v\,{\rm mod}(n_j)}W_{v'}\ \ .
\end{equation}
Also we get the following identifications of real vector bundles over $F$
(cf. \cite[(4.11)]{MR1870666}),
\begin{equation}
\begin{split}\label{s27}
TX(n_j)\big|_{F}&=TY\oplus\bigoplus_{v>0,\atop v\equiv0\,{\rm mod}\,(n_j)}N_v\ ,\quad
N(n_j)_{\frac{n_j}{2}}^{\mathbb{R}}\Big|_{F}=\bigoplus_
{v>0,\atop v\equiv\frac{n_j}{2}\,{\rm mod}\,(n_j)}N_v\ \ ,
\\
V(n_j)_0^{\mathbb{R}}\Big|_{F}&=V_0^{\mathbb{R}}\oplus\bigoplus_{v>0,\atop v\equiv0\,{\rm mod}(n_j)}V_v\ ,\quad
V(n_j)_{\frac{n_j}{2}}^{\mathbb{R}}\Big|_{F}=\bigoplus_{v>0,\atop v\equiv\frac{n_j}{2}\,{\rm mod}(n_j)}V_v\ \ .
\end{split}
\end{equation}
Moreover, we have the identifications of complex vector bundles over $F$,
\begin{equation}\label{eg51}
\begin{split}
TX(n_j)\big|_{F}\otimes_{\mathbb{R}}\mathbb{C}&=TY\otimes_{\mathbb{R}}\mathbb{C}\oplus
\bigoplus_{v>0,\,v\equiv0\,{\rm mod}(n_j)}(N_v\oplus\overline{N}_v)\ ,
\\
V(n_j)_0\big|_{F}&=V_0^{\mathbb{R}}\otimes_{\mathbb{R}}\mathbb{C}\oplus
\bigoplus_{v>0,\,v\equiv0\,{\rm mod}(n_j)}(V_v\oplus\overline{V}_v)\ .
\end{split}
\end{equation}

As $(p_j,n_j)=1$, we know that, for $v\in\mathbb{Z}$, $\tfrac{p_j}{n_j}v\in\mathbb{Z}$ if and only if $\tfrac{v}{n_j}\in\mathbb{Z}$.
Also, $\tfrac{p_j}{n_j}v\in\mathbb{Z}+\tfrac{1}{2}$ if and only if $\tfrac{v}{n_j}\in\mathbb{Z}+\tfrac{1}{2}$. Remark if
$v\equiv -v'\ {\rm mod}(n_j)$, then $\{n\,|\,0<n\in\mathbb{Z}+\tfrac{p_j}{n_j}v\}=
\{n\,|\,0<n\in\mathbb{Z}-\tfrac{p_j}{n_j}v'\}$. Using the identifications
(\ref{uy2}), (\ref{uy3}) and (\ref{eg51}), we can rewrite $\mathcal{F}(\beta_j)$, $F_V^1(\beta_j)$,
$F_V^2(\beta_j)$, $Q^1_W(\beta_j)$ and $Q^2_W(\beta_j)$ over $F$ defined over in (\ref{s18}) as follows (compare with \cite[(4.7)]{MR1870666}),

\begin{equation}
\begin{split}\label{uy8}
\mathcal{F}(\beta_j)&=\bigotimes_{0<n\in\mathbb{Z}}{\rm Sym}\,\bigl(TX(n_j)_n\bigr)\otimes\bigotimes_{0<v<n_j/2}{\rm Sym}\,
\Bigl(\bigoplus_{0<n\in\mathbb{Z}+\frac{p_j}{n_j}v}N(n_j)_{v,n}\Bigr.
\\
&\hspace{30pt}
\Bigl.\oplus\bigoplus_{0<n\in\mathbb{Z}-\frac{p_j}{n_j}v}
\overline{N(n_j)}_{v,n}\Bigr)\oplus\bigoplus_{0<n\in\mathbb{Z}+\frac{1}{2}}{\rm Sym}\,\bigl(N(n_j)_{n_j/2,n}\bigr)\ ,
\end{split}
\end{equation}
\begin{equation}
\begin{split}\label{uy9}
F_V^1(\beta_j)&=\Lambda\Bigl(\ \bigoplus_{0<n\in\mathbb{Z}}V(n_j)_{0,n}\oplus\bigoplus_{0<v<n_j/2}
\Bigl(\ \bigoplus_{0<n\in\mathbb{Z}+\frac{p_j}{n_j}v}V(n_j)_{v,n}\Bigr.\Bigr.
\\
&\hspace{30pt}
\Bigl.\Bigl.\oplus\bigoplus_{0<n\in\mathbb{Z}-\frac{p_j}{n_j}v}
\overline{V(n_j)}_{v,n}\Bigr)\bigoplus_{0<n\in\mathbb{Z}+\frac{1}{2}}V(n_j)_{n_j/2,n}\Bigr)\ ,
\end{split}
\end{equation}
\begin{equation}
\begin{split}\label{uy10}
F_V^2(\beta_j)&=\Lambda\Bigl(\ \bigoplus_{0<n\in\mathbb{Z}}V(n_j)_{n_j/2,n}\oplus\bigoplus_{0<v<n_j/2}
\Bigl(\ \bigoplus_{0<n\in\mathbb{Z}+\frac{p_j}{n_j}v+\frac{1}{2}}V(n_j)_{v,n}\Bigr.\Bigr.
\\
&\hspace{30pt}
\Bigl.\Bigl.\oplus\bigoplus_{0<n\in\mathbb{Z}-\frac{p_j}{n_j}v+\frac{1}{2}}\overline{V(n_j)}_{v,n}\Bigr)
\bigoplus_{0<n\in\mathbb{Z}+\frac{1}{2}}V(n_j)_{0,n}\Bigr)\ ,
\end{split}
\end{equation}
\begin{equation}
Q_{W}^1(\beta_j)=\Lambda\Bigl(\ \bigoplus_{0\leq v<n_j}\Bigl(\bigoplus_
{0<n\in\mathbb{Z}+\frac{p_j}{n_j}v}W(n_j)_{v,n}\oplus\bigoplus_{0\leq n\in\mathbb{Z}-\frac{p_j}{n_j}v}
\overline{W(n_j)}_{v,n}\Bigr)\,\Bigr)\ ,
\end{equation}
\begin{equation}
\begin{split}\label{qw2}
Q_{W}^2(\beta_j)=\Lambda\Big(\ \bigoplus_{0\leq v<n_j}\Big(&\bigoplus_
{0<n\in\mathbb{Z}+\frac{p_j}{n_j}v+{\frac{1}{2}}}W(n_j)_{v,n}
\\
&\hspace{40pt}\oplus\bigoplus_{0<n\in\mathbb{Z}-\frac{p_j}{n_j}v+{\frac{1}{2}}}
\overline{W(n_j)}_{v,n}\Big)\,\Big)\ .
\end{split}
\end{equation}
We indicate here that $\mathcal{F}(\beta_j)$, $F_V^1(\beta_j)$,
$F_V^2(\beta_j)$, $Q_{W}^1(\beta_j)$ and $Q_{W}^2(\beta_j)$ in (\ref{s18})
are the restriction of the corresponding $\mathbb{Z}/k$ vector bundles in the right
hand side of \eqref{uy8}-\eqref{qw2} over $M(n_j)$, which will
still be denoted as $\mathcal{F}(\beta_j)$, $F_V^1(\beta_j)$, $F_V^2(\beta_j)$, $Q_{W}^1(\beta_j)$ and $Q_{W}^2(\beta_j)$.
Write
\begin{equation}\label{s29}
Q_W(\beta_j)=Q^1_W(\beta_j)\otimes Q^2_W(\beta_j)\otimes Q^2_W(\beta_j),
\end{equation}
which is a $\mathbb{Z}/k$ vector bundles over  $M(n_j)$.

We now define the Spin$^c$ Dirac operators on $M(n_j)$.
The following Lemma follows from the proof of  \cite[Lemmas 11.3 and 11.4]{MR954493}.

\begin{lem}\label{eg54}
{\rm (Compare with \cite[Lemma 4.2]{MR1870666})} Assume that {\rm \eqref{hypothesis}} holds. Let
\begin{equation}\label{eg55}
\begin{split}
L(n_j)&=\bigotimes_{0<v<n_j/2}\left(\det (N(n_j)_v)\otimes\det(\overline{V(n_j)}_v)\right.
\\&\hspace{60pt}\left.
\otimes\Bigl(\det(\overline{W(n_j)}_v)
\otimes\det(W(n_j)_{n_j-v})\Bigr)^{3}\,\right)^{(r(n_j)+1)v}
\end{split}
\end{equation}
be the complex line bundle over $M(n_j)$. Then we have
\begin{enumerate}[{\rm (i)}]
\item $L(n_j)$ has an $n_j^{\rm th}$ root over $M(n_j)$.

\item Let\ \ $U_1=TX(n_j)\oplus V(n_j)_0^{\mathbb{R}}\oplus W(n_j)_{{n_j}/{2}}^{\mathbb{R}}\oplus W(n_j)_{{n_j}/{2}}^{\mathbb{R}}$\ ,
\begin{equation*}
\begin{split}
L_1&=K_X\otimes\bigotimes_{0<v<n_j/2}\Bigl(\det (N(n_j)_v)\otimes\det(\overline{V(n_j)}_v)\Bigr)
\\
&\hspace{90pt}
\otimes\Bigl(\det(W(n_j)_{n_j/2})\Bigr)^3
\otimes L(n_j)^{\frac{r(n_j)}{n_j}}\ .
\end{split}
\end{equation*}
Let\ \ $U_2=TX(n_j)\oplus V(n_j)_{n_j/2}^{\mathbb{R}}\oplus W(n_j)_{{n_j}/{2}}^{\mathbb{R}}\oplus W(n_j)_{{n_j}/{2}}^{\mathbb{R}}$\ ,
\begin{equation*}
L_2=K_X\otimes\bigotimes_{0<v<n_j/2}\Bigl(\det (N(n_j)_v)\Bigr)
\otimes\Bigl(\det(W(n_j)_{n_j/2})\Bigr)^3\otimes L(n_j)^{\frac{r(n_j)}{n_j}}.
\end{equation*}
Then
$U_1$ (resp. $U_2$) has a Spin$^c$ structure defined by $L_1$ (resp. $L_2$).
\end{enumerate}
\end{lem}

Remark that in order to define an $S^1$ (resp. $G_y$) action on $L(n_j)^{r(n_j)/n_j}$, we must replace the $S^1$
(resp.  $G_y$) action by its $n_j$-fold action. Here by abusing notation,
we still say an $S^1$ (resp. $G_y$) action without causing any confusion.

Let $S(U_1,L_1)$ (resp. $S(U_2,L_2)$) be the fundamental complex spinor bundle for
$(U_1,L_1)$ (resp. $(U_2,L_2)$) (cf. \cite[Appendix D]{MR1031992}). There are two $\mathbb{Z}_2$-gradings
on $S(U_1,L_1)$ (resp. $S(U_2,L_2)$). The first grading, which we denote by $\tau_s$, is induced by
the involutions on $S(U_1,L_1)$ and $S(U_2,L_2)$ determined by
$TX(n_j)\oplus W(n_j)_{{n_j}/{2}}^{\mathbb{R}}$ as in \eqref{s5}.
The second grading, which we denote by $\tau_e$, is induced by
the involution on $S(U_1,L_1)$ (resp. $S(U_2,L_2)$) determined by
$TX(n_j)\oplus V(n_j)_0^{\mathbb{R}}\oplus W(n_j)_{{n_j}/{2}}^{\mathbb{R}}$
(resp. $U_2=TX(n_j)\oplus V(n_j)_{n_j/2}^{\mathbb{R}}\oplus W(n_j)_{{n_j}/{2}}^{\mathbb{R}}$) as in \eqref{s5}.

In what follows, by $D^{X(n_j)}$ we mean the $S^1$-equivariant Spin$^c$ Dirac operator on $S(U_1,L_1)$ or $S(U_2,L_2)$
over $M(n_j)$.

Corresponding to \eqref{s21}, by (\ref{uy2}) and \eqref{uy3}, we define $S(U_1,L_1)'$ (resp. $S(U_2,L_2)'$)
equipped with its involutions $\tau'_s$ and $\tau'_e$ as follows (compare with \cite[(4.16)]{MR1870666}),
\begin{align}
\label{uy5}
&\Bigl(S(U_1,L_1)',\tau'_s/\tau'_e\Bigr)=\Big(S\Big(TY\oplus V_0^{\mathbb{R}},L_1\otimes
\bigotimes_{v>0,\atop v\equiv0\,{\rm mod}\,(n_j)}
(\det N_v\otimes\det V_v)^{-1}
\notag\\
&\hspace{60pt}\otimes\bigotimes_{v\equiv\frac{n_j}{2}\,{\rm mod}\,(n_j)}(\det W_v)^{-2}\Big),\tau'_s/\tau'_e\Big)
\otimes\bigotimes_{v>0,\atop v\equiv0\,{\rm mod}\,(n_j)}\Lambda_{\pm 1} (V_v)
\notag\\
&\hspace{40pt}\otimes\bigotimes_{v\equiv \frac{n_j}{2}\,{\rm mod}\,(n_j)}\Lambda_{-1}(W_v)
\otimes\bigotimes_{v\equiv \frac{n_j}{2}\,{\rm mod}\,(n_j)}\Lambda(W_v)\ ,
\end{align}
\begin{align}
\label{uy6}
&\Big(S(U_2,L_2)',\tau'_s/\tau'_e\Bigr)=S\Bigl(TY,L_2
\otimes\bigotimes_{v>0,\atop v\equiv0\ {\rm mod}(n_j)}(\det N_v)^{-1}
\notag\\
&\hspace{60pt}\otimes\bigotimes_{v>0,\atop v\equiv \frac{n_j}{2}\,{\rm mod}\,(n_j)}(\det V_v)^{-1}
\otimes\bigotimes_{v\equiv\frac{n_j}{2}\,{\rm mod}\,(n_j)}(\det W_v)^{-2}\Big)
\\
&\hspace{20pt}\otimes\bigotimes_{v>0,\atop v\equiv\frac{n_j}{2}\,{\rm mod}\,(n_j)}\Lambda_{\pm 1} (V_v)
\otimes\bigotimes_{v\equiv \frac{n_j}{2}\,{\rm mod}\,(n_j)}\Lambda_{-1}(W_v)
\otimes\bigotimes_{v\equiv \frac{n_j}{2}\,{\rm mod}\,(n_j)}\Lambda(W_v)\ .
\notag
\end{align}
Then by \eqref{s21}, for $i=1, 2$, we have the following isomorphisms of
Clifford modules over $F$ preserving the $\mathbb{Z}_2$-gradings (compare with \cite[(4.17)]{MR1870666}),
\begin{equation}\label{eg58}
\Bigl(S(U_i,L_i),\tau_s/\tau_e\Bigr)\Big|_{F}\simeq \Bigl(S(U_i,L_i)',\tau'_s/\tau'_e\Bigr)
\otimes\bigotimes_{v>0,\,v\equiv0\,{\rm mod}\,(n_j)}\Lambda_{-1} (N_v).
\end{equation}

As in \cite[pp. 952]{MR1870666}, we introduce formally the following complex line bundles over $F$,
\begin{align}\label{s23}
L_1'&=\Big(L_1^{-1}\otimes\bigotimes_{ v>0,\atop v\equiv0\,{\rm mod}\,(n_j)}
(\det N_v\otimes\det V_v)\otimes
\bigotimes_{ v\equiv \frac{n_j}{2}\,{\rm mod}\,(n_j)}(\det W_{v})^2
\notag\\
&\qquad\otimes\bigotimes_{v>0}\,(\det N_v\otimes\det V_v)^{-1}\otimes K_X\Big)^{\frac{1}{2}},
\end{align}
\begin{align}\label{s24}
L_2'&=\Big(L_2^{-1}
\otimes\bigotimes_{v>0,\atop v\equiv0\,{\rm mod}\,(n_j)}\det N_v
\otimes\bigotimes_{v>0,\atop v\equiv\frac{n_j}{2}\,{\rm mod}\,(n_j)}\det V_v
\notag\\
&\qquad\otimes\bigotimes_{v\equiv \frac{n_j}{2}\,{\rm mod}\,(n_j)}(\det W_{v})^2
\otimes\bigotimes_{v>0}\,(\det N_v)^{-1}\otimes K_X\Big)^{\frac{1}{2}}.
\end{align}
In fact, from \eqref{s21}, Lemma \ref{eg54} and the assumption that $V$ is spin,
one verifies easily that $c_1(L_i'^2)=0 \mod (2)$ for $i=1,2$, which implies that $L_1'$ and $L_2'$ are
well-defined complex line bundles over $F$.

Then by \cite[(3.14)]{MR1870666}, \eqref{uy5}-\eqref{s24}
and the definitions of $L_1$, $L_2$, we get the
following identifications of Clifford modules over $F$ (compare with \cite[(4.19)]{MR1870666}),
\begin{align}\label{eg60}
&\Bigl(S(U_1,L_1)'\otimes L_1',(\tau'_s/\tau'_e)\otimes {\rm id}\Bigr)=S\bigl(TY,K_X\otimes_{v>0}(\det N_v)^{-1}\bigr)
\notag\\
&\hspace{60pt}\otimes\Bigl(S\bigl(V_0^{\mathbb{R}},\otimes_{v>0}(\det V_v)^{-1}\bigr),{\rm id}/\tau\Bigr)
\otimes\bigotimes_{v>0,\atop v\equiv0\,{\rm mod}\,(n_j)}\Lambda_{\pm 1}(V_v)
\notag\\
&\hspace{40pt}\otimes\bigotimes_{v>0, \atop v\equiv \frac{n_j}{2}\,{\rm mod}\,(n_j)}\Lambda_{-1}(W_v)
\otimes\bigotimes_{v<0, \atop v\equiv \frac{n_j}{2}\,{\rm mod}\,(n_j)}\Lambda_{-1}(\overline{W}_v)
\\
&\otimes\bigotimes_{v>0, \atop v\equiv \frac{n_j}{2}\,{\rm mod}\,(n_j)}\Lambda(W_v)
\otimes\bigotimes_{v<0, \atop v\equiv \frac{n_j}{2}\,{\rm mod}\,(n_j)}\Lambda(\overline{W}_v)
\otimes
\bigotimes_{v<0,\atop v\equiv \frac{n_j}{2}{\rm mod}\,(n_j)}(\det W_v)^2\ ,
\notag
\end{align}
\begin{align}
\label{uy30}
&\Bigl(S(U_2,L_2)'\otimes L_2',(\tau'_s/\tau'_e)\otimes {\rm id}\Bigr)
=S\bigl(TY,K_X\otimes_{v>0}(\det N_v)^{-1}\bigr)
\notag\\
&\hspace{20pt}\otimes\bigotimes_{v>0,\atop v\equiv \frac{n_j}{2}
\,{\rm mod}\,(n_j)}\Lambda_{\pm 1}(V_v)
\otimes\bigotimes_{v>0, \atop v\equiv \frac{n_j}{2}\,{\rm mod}\,(n_j)}\Lambda_{-1}(W_v)
\otimes\bigotimes_{v<0, \atop v\equiv \frac{n_j}{2}\,{\rm mod}\,(n_j)}\Lambda_{-1}(\overline{W}_v)
\notag\\
&\hspace{10pt}\otimes\bigotimes_{v>0, \atop v\equiv \frac{n_j}{2}\,{\rm mod}\,(n_j)}\Lambda(W_v)
\otimes\bigotimes_{v<0, \atop v\equiv \frac{n_j}{2}\,{\rm mod}\,(n_j)}\Lambda(\overline{W}_v)
\otimes\bigotimes_{v<0,\atop v\equiv \frac{n_j}{2}{\rm mod}\,(n_j)}(\det W_v)^2\ .
\end{align}

Now we compare the $\mathbb{Z}_2$-gradings in (\ref{eg60}) and \eqref{uy30}.
Set (compare with \cite[(4.20)]{MR1870666})
\begin{equation}\label{uy7}
\begin{split}
\Delta(n_j,N)&=\sum_{\frac{n_j}{2}<v'<n_j}\sum_{0<v,\,v\equiv v'\,{\rm mod}\,(n_j)}
\dim N_v+o\Bigl(N(n_j)_{\frac{n_j}{2}}^{\mathbb{R}}\Bigr)\ ,
\\
\Delta(n_j,V)&=\sum_{\frac{n_j}{2}<v'<n_j}\sum_{0<v,\,v\equiv v'\,{\rm mod}\,(n_j)}
\dim V_v+o\Bigl(V(n_j)_{\frac{n_j}{2}}^{\mathbb{R}}\Bigr)\ ,
\\
\Delta(n_j,W)&=\sum_{v<0,\,v\equiv \frac{n_j}{2}\,{\rm mod}\,(n_j)}\dim W_v\ ,
\end{split}
\end{equation}
with $o(N(n_j)_{\frac{n_j}{2}}^{\mathbb{R}})$ (resp. $o(V(n_j)_{\frac{n_j}{2}}^{\mathbb{R}})$ equals 0 or 1,
depending on whether the given orientation on $N(n_j)_{\tfrac{n_j}{2}}^{\mathbb{R}}$
(resp. $V(n_j)_{\tfrac{n_j}{2}}^{\mathbb{R}}$) agrees or disagrees with the complex
orientation of $\oplus_{v>0,\,v\equiv \frac{n_j}{2}\,{\rm mod}\,(n_j)} N_v$
(resp. $\oplus_{v>0,\,v\equiv \frac{n_j}{2}\,{\rm mod}\,(n_j)} V_v$).

As explained in \cite[pp. 166]{MR2016198}, for the $\mathbb{Z}_2$-gradings induced by $\tau_s$,
the differences of the $\mathbb{Z}_2$-gradings of (\ref{eg60}) and \eqref{uy30} are both $(-1)^{\Delta(n_j,N)+\Delta(n_j,W)}$;
for the $\mathbb{Z}_2$-gradings induced by $\tau_e$,
the difference of the $\mathbb{Z}_2$-gradings of (\ref{eg60}) (resp. \eqref{uy30}) is
$(-1)^{\Delta(n_j,N)+\Delta(n_j,V)+\Delta(n_j,W)}$ (resp. $\,(-1)^{\Delta(n_j,N)+o(V(n_j)_{n_j/2}^{\mathbb{R}})+\Delta(n_j,W)}\,$).

\begin{lem}\label{eg61}
{\rm (Compare with \cite[Lemma 4.3]{MR1870666})} Let us write
\begin{equation*}
\begin{split}
L(\beta_j)_1&=L_1'
\otimes\bigotimes_{v>0}(\det N_v)^{\bigl[\frac{p_j}{n_j}v\bigr]+(p-1)v+1}
\otimes\bigotimes_{v>0}(\det \overline{V}_v)^{\bigl[\frac{p_j}{n_j}v\bigr]+(p-1)v}
\\
&\otimes\bigotimes_{v>0,\atop v\equiv0\,{\rm mod}(n_j)}(\det N_v)^{-1}
\otimes\bigotimes_{v<0}(\det W_v)^{\bigl[-\frac{p_j}{n_j}v\bigr]+2\bigl[-\frac{p_j}{n_j}v+\frac{1}{2}\bigr]-3(p-1)v}
\\
&\otimes\bigotimes_{v>0}\,(\det \overline{W}_v)^{\bigl[\frac{p_j}{n_j}v\bigr]+2\bigl[\frac{p_j}{n_j}v+\frac{1}{2}\bigr]+3(p-1)v+1}
\\
&\otimes\bigotimes_{v>0,\atop v\equiv0\,{\rm mod}(n_j)}\det {W}_v
\otimes\bigotimes_{v<0,\atop v\equiv\frac{n_j}{2}\,{\rm mod}(n_j)}(\det \overline{W}_v)^2\  ,
\end{split}
\end{equation*}
\begin{equation*}
\begin{split}
L(\beta_j)_2&=L_2'
\otimes\bigotimes_{v>0}(\det N_v)^{\bigl[\frac{p_j}{n_j}v\bigr]+(p-1)v+1}
\otimes\bigotimes_{v>0}(\det \overline{V}_v)^{\bigl[\frac{p_j}{n_j}v+\frac{1}{2}\bigr]+(p-1)v}
\\&
\otimes\bigotimes_{v>0,\atop v\equiv0\,{\rm mod}(n_j)}(\det N_v)^{-1}
\otimes\bigotimes_{v<0}(\det W_v)^{\bigl[-\frac{p_j}{n_j}v\bigr]+2\bigl[-\frac{p_j}{n_j}v+\frac{1}{2}\bigr]-3(p-1)v}
\\
&\otimes\bigotimes_{v>0}\,(\det \overline{W}_v)^{\bigl[\frac{p_j}{n_j}v\bigr]+2\bigl[\frac{p_j}{n_j}v+\frac{1}{2}\bigr]+3(p-1)v+1}
\\
&\otimes\bigotimes_{v>0,\atop v\equiv0\,{\rm mod}(n_j)}\det {W}_v
\otimes\bigotimes_{v<0,\atop v\equiv\frac{n_j}{2}\,{\rm mod}(n_j)}(\det \overline{W}_v)^2\ .
\end{split}
\end{equation*}
Then $L(\beta_j)_1$ and $L(\beta_j)_2$ can be extended naturally to $G_y\times S^1$-equivariant complex
line bundles over $M(n_j)$ which we will still denote by $L(\beta_j)_1$ and $L(\beta_j)_2$  respectively.
\end{lem}
{\bf Proof}\hspace{4mm} We introduce the following line bundle over $M(n_j)$,
\begin{equation}
\begin{split}
L^\omega(\beta_j)&=\bigotimes_{0<v<n_j/2}\Big(\det (N(n_j)_v)\otimes\det(\overline{V(n_j)}_v)
\\&\hspace{40pt}
\otimes\bigl(\det(\overline{W(n_j)}_v)
\otimes\det(W(n_j)_{n_j-v})\bigr)^{3}\,\Big)^{-\omega(v)-r(n_j)v}\ .
\end{split}
\end{equation}
where as in \cite[(4.35)]{MR2016198} we write $[\tfrac{p_j}{n_j}v]=\tfrac{p_j}{n_j}v-\tfrac{\omega(v)}{n_j}$\ .

As in \cite[(4.38)]{MR2016198} and \cite[(4.28)]{MR1870666}, Lemma \ref{eg54} implies that $L^\omega(\beta_j)^{1/n_j}$ is
well-defined over $M(n_j)$. Direct calculation shows that
\begin{equation*}
\begin{split}
L(\beta_j)_1&=L^{-(p-1)-\frac{p_j}{n_j}}\otimes L^\omega(\beta_j)^{\frac{1}{n_j}}
\otimes \bigotimes_{0<v<\frac{n_j}{2}}\det(\overline{W(n_j)}_v)
\otimes\bigl(\det(\overline{W(n_j)}_{\frac{n_j}{2}})\bigr)^2
\\
&\otimes\bigotimes_{1\leq m\leq \frac{p_j}{2}}\bigotimes_{m-\frac{1}{2}<\frac{p_j}{n_j}v<m}
\bigl(\det(\overline{W(n_j)}_v)\otimes\det(W(n_j)_{n_j-v})\bigr)^2\ ,
\end{split}
\end{equation*}
\begin{equation*}
\begin{split}
&L(\beta_j)_2=L^{-(p-1)-\frac{p_j}{n_j}}\otimes L^\omega(\beta_j)^{\frac{1}{n_j}}
\otimes \bigotimes_{0<v<\frac{n_j}{2}}\det(\overline{W(n_j)}_v)
\otimes\bigl(\det(\overline{W(n_j)}_{\frac{n_j}{2}})\bigr)^2
\\
&\otimes\bigotimes_{1\leq m\leq \frac{p_j}{2}}\bigotimes_{m-\frac{1}{2}<\frac{p_j}{n_j}v<m}
\left(\bigl(\det(\overline{W(n_j)}_v)\otimes\det(W(n_j)_{n_j-v})\bigr)^2\otimes\det(\overline{V(n_j)}_v)\right).
\end{split}
\end{equation*}

The proof of Lemma \ref{eg61} is completed.

\

To simplify the notations, we introduce the following locally constant functions on $F$
(compare with \cite[(4.45)]{MR2016198}, \cite[(4.30)]{MR1870666}),
\begin{equation}
\begin{split}\label{s41}
\varepsilon^1_W&=-\frac{1}{2}\sum_{v>0}(\dim W_v)\cdot\Bigl(\bigl(\bigl[\tfrac{p_j}{n_j}v\bigr]+(p-1)v\bigr)
\bigl(\bigl[\tfrac{p_j}{n_j}v\bigr]+(p-1)v+1\bigr)\Bigr.
\\
&\hspace{60pt}
-\Bigl.(\tfrac{p_j}{n_j}v+(p-1)v)\bigl(2\bigl(\bigl[\tfrac{p_j}{n_j}v\bigr]+(p-1)v\bigr)+1\bigr)\Bigr)
\\
&-\frac{1}{2}\sum_{v<0}(\dim W_v)\cdot\Bigl(\bigl(\bigl[-\tfrac{p_j}{n_j}v\bigr]-(p-1)v\bigr)
\bigl(\bigl[-\tfrac{p_j}{n_j}v\bigr]-(p-1)v+1\bigr)\Bigr.
\\
&\hspace{60pt}
+\Bigl.\bigl(\tfrac{p_j}{n_j}v+(p-1)v\bigr)\bigl(2\bigl(\bigl[-\tfrac{p_j}{n_j}v\bigr]-(p-1)v\bigr)+1\bigr)\Bigr)\ ,
\end{split}
\end{equation}
\begin{equation}
\begin{split}\label{s40}
\varepsilon^2_W&=-\frac{1}{2}\sum_{v>0}(\dim W_v)\cdot\Bigl(
\bigl(\bigl[\tfrac{p_j}{n_j}v+\tfrac{1}{2}\bigr]+(p-1)v\bigr)^2
\Bigr.
\\
&\hspace{60pt}-\Bigl.
2\,(\tfrac{p_j}{n_j}v+(p-1)v)\bigl(\bigl[\tfrac{p_j}{n_j}v+\tfrac{1}{2}\bigr]+(p-1)v\bigr)
\Bigr)
\\
&-\frac{1}{2}\sum_{v<0}(\dim W_v)\cdot\Bigl(
\bigl(\bigl[-\tfrac{p_j}{n_j}v+\tfrac{1}{2}\bigr]-(p-1)v\bigr)^2
\Bigr.
\\
&\hspace{60pt}+\Bigl.
2\,(\tfrac{p_j}{n_j}v+(p-1)v)\bigl(\bigl[-\tfrac{p_j}{n_j}v+\tfrac{1}{2}\bigr]-(p-1)v\bigr)
\Bigr)\ ,
\end{split}
\end{equation}

\begin{align}
\varepsilon_1&=\frac{1}{2}\sum_{v>0}(\dim N_v-\dim V_v)\Big(\bigl(\bigl[\tfrac{p_j}{n_j}v\bigr]+(p-1)v\bigr)
\bigl(\bigl[\tfrac{p_j}{n_j}v\bigr]+(p-1)v+1\bigr)\Bigr.
\notag\\
&\hspace{60pt}
\Bigl.-\bigl(\tfrac{p_j}{n_j}v+(p-1)v\bigr)\bigl(2\bigl(\bigl[\tfrac{p_j}{n_j}v\bigr]+(p-1)v\bigr)+1\bigr)\Bigr)\ ,
\end{align}
\begin{equation}
\begin{split}
\varepsilon_2&=\frac{1}{2}\sum_{v>0}(\dim N_v)\cdot\Bigl(\bigl(\bigl[\tfrac{p_j}{n_j}v\bigr]+(p-1)v\bigr)
\bigl(\bigl[\tfrac{p_j}{n_j}v\bigr]+(p-1)v+1\bigr)\Bigr.
\\
&\hspace{70pt}
\Bigl.-(\tfrac{p_j}{n_j}v+(p-1)v)\bigl(2\bigl(\bigl[\tfrac{p_j}{n_j}v\bigr]+(p-1)v\bigr)+1\bigr)\Bigr)
\\
&\quad-\frac{1}{2}\sum_{v>0}(\dim V_v)\cdot\Bigl(\bigl(\bigl[\tfrac{p_j}{n_j}+\tfrac{1}{2}\bigr]+(p-1)v\bigl)^2\Bigr.
\\
&\hspace{70pt}
\Bigl.-2\bigl(\tfrac{p_j}{n_j}v+(p-1)v\bigr)\bigl(\bigl[\tfrac{p_j}{n_j}+\tfrac{1}{2}\bigr]+(p-1)v\bigr)\Bigr)\ .
\end{split}
\end{equation}
As in \cite[(2.23)]{MR1870666}, for $0\leq j\leq J_0$, we set
\begin{align}\label{eg64}
e(p,\beta_j,N)&=\frac{1}{2}\sum_{v>0}(\dim N_v)\cdot\bigl(\bigl[\tfrac{p_j}{n_j}v\bigr]+(p-1)v\bigr)
\bigl(\bigl[\tfrac{p_j}{n_j}v\bigl]+(p-1)v+1\bigl)\ ,
\notag\\
d'(p,\beta_j,N)&=\sum_{v>0}(\dim N_v)\cdot\bigl(\bigl[\tfrac{p_j}{n_j}v\bigr]+(p-1)v\bigr)\ .
\end{align}
Then $e(p,\beta_j,N)$ and $d'(p,\beta_j,N)$ are locally constant functions on $F$. In particular, we have
\begin{equation}\label{eg65}
\begin{split}
e(p,\beta_0,N)&=\frac{1}{2}(p-1)^2e(N)+\frac{1}{2}(p-1)d'(N)\ ,
\\
e(p,\beta_{J_0},N)&=\frac{1}{2}p^2e(N)+\frac{1}{2}p\,d'(N)\ ,
\\
d'(p,\beta_{J_0},N)&=d'(p+1,\beta_{0},N)=p\,d'(N)\ .
\end{split}
\end{equation}

\begin{prop}{\rm (Compare with \cite[Proposition 4.2]{MR1870666})} \label{s25}
For $i=1,2$, the $G_y$-equivariant isomorphisms of complex vector bundles over $F$ induced by Proposition \ref{eg44} and
\eqref{eg60}-\eqref{uy30},
\begin{equation*}
\begin{split}
&r_{i1}:\,S(TY,K_X\otimes_{v>0}(\det N_v)^{-1})
\otimes(K_W\otimes K^{-1}_X)^{1/2}
\\
&\hspace{60pt}\otimes \mathcal{F}_{p,j-1}(X)\otimes F^i_V\otimes Q(W)
\\
&\qquad\longrightarrow S(U_i,L_i)'
\otimes(K_W\otimes K^{-1}_X)^{1/2}\otimes \mathcal{F}(\beta_j)\otimes F_V^i(\beta_j)
\\&\hspace{60pt}
\otimes Q_W(\beta_j)\otimes L(\beta_j)_i\otimes\bigotimes_{v>0,\atop v\equiv0\,{\rm mod}\,(n_j)}{\rm Sym}\,(\overline{N}_{v,0})\ ,
\end{split}
\end{equation*}
\begin{equation*}
\begin{split}
&r_{i2}:\,S(TY,K_X\otimes_{v>0}(\det N_v)^{-1})
\otimes(K_W\otimes K^{-1}_X)^{1/2}
\\
&\hspace{60pt}\otimes \mathcal{F}_{p,j}(X)\otimes F^i_V\otimes Q(W)
\\
&\qquad\longrightarrow S(U_i,L_i)'
\otimes(K_W\otimes K^{-1}_X)^{1/2}\otimes \mathcal{F}(\beta_j)\otimes F_V^i(\beta_j)
\\&\hspace{60pt}
\otimes Q_W(\beta_j)\otimes L(\beta_j)_i\otimes
\bigotimes_{v>0,\atop v\equiv0\,{\rm mod}\,(n_j)}\Bigl({\rm Sym}\,({N}_{v,0})\otimes \det N_v\Bigr) ,
\end{split}
\end{equation*}
have the following properties:
\begin{enumerate}[{\rm (i)}]
\item for $i=1,2$, $\gamma=1,2$,
\begin{equation}\label{eg68}
\begin{split}
r_{i\gamma}^{-1}\cdot \textbf{J}_H\cdot r_{i\gamma}&=\textbf{J}_H\ ,
\\
r_{i\gamma}^{-1}\cdot P\cdot r_{i\gamma}&=P+\Bigl(\,\frac{p_j}{n_j}+(p-1)\Bigr)\textbf{J}_H+\varepsilon_{i\gamma}\ ,
\end{split}
\end{equation}
where $\varepsilon_{i\gamma}$ are given by
\begin{equation}
\begin{split}
\varepsilon_{i1}&=\varepsilon_i+\varepsilon^1_W+2\varepsilon^2_W-e(p,\beta_{j-1},N)\ ,
\\
\varepsilon_{i2}&=\varepsilon_i+\varepsilon^1_W+2\varepsilon^2_W-e(p,\beta_{j},N)\ .
\end{split}
\end{equation}

\item for $i=1,2$, $\gamma=1,2$,
\begin{equation}\label{eg71}
\begin{split}
r_{i\gamma}^{-1}\tau_{e}r_{i\gamma}&=(-1)^{\mu_i}\tau_{e}\ ,\quad r_{i\gamma}^{-1}\tau_{s}r_{i\gamma}
=(-1)^{\mu_3}\tau_{s}\ ,
\\
r_{i\gamma}^{-1}\tau_{1}r_{i\gamma}&=(-1)^{\mu_4}\tau_{1}\ ,
\end{split}
\end{equation}
where $\mu_{i}$ are given by
\begin{equation*}
\begin{split}
\mu_1&=-\sum_{v>0}(\dim V_v)\bigl[\tfrac{p_j}{n_j}v\bigr]+\Delta(n_j,N)+\Delta(n_j,V)+\Delta(n_j,W) \mod (2),
\\
\mu_2&=-\sum_{v>0}(\dim V_v)\cdot\bigl[\tfrac{p_j}{n_j}v+\tfrac{1}{2}\bigr]+\Delta(n_j,N)
\\
&\hspace{70pt}+o\bigl(V(n_j)_{\tfrac{n_j}{2}}^{\mathbb{R}}\bigr)+\Delta(n_j,W) \mod (2),
\\
\mu_3&=\Delta(n_j,N)+\Delta(n_j,W) \mod (2),
\\
\mu_4&=\sum_{v>0}(\dim W_v)\cdot\bigl(\bigl[\tfrac{p_j}{n_j}v+\tfrac{1}{2}\bigr]+(p-1)v\bigr)
\\
&\hspace{40pt}+\sum_{v<0}(\dim W_v)\cdot\bigl(\bigl[-\tfrac{p_j}{n_j}v+\tfrac{1}{2}\bigr]-(p-1)v\bigr) \mod (2).
\end{split}
\end{equation*}
\end{enumerate}
\end{prop}
{\bf Proof}\hspace{4mm}By the proof of  \cite[Proposition 4.2]{MR1870666}, we need to compute the action of
$r_*^{-1}\cdot P\cdot r_*$ on $\bigotimes_{0<n\in\mathbb{Z}+{1\over 2},\,v>0,\atop n-(p-1)v-\frac{p_j}{n_j}v\leq 0}
\Lambda^{i_n}(\overline{W}_{v,n})
\otimes
\bigotimes_{0<n\in\mathbb{Z}+{1\over 2},\,v<0,\atop n+(p-1)v+\frac{p_j}{n_j}v\leq 0}
\Lambda^{i'_n}({W}_{v,n})$.
In fact, by \eqref{s17}, as in \eqref{s13}, we get
\begin{equation}
\begin{split}\label{s42}
r_*^{-1}\cdot P\cdot r_*&=
\sum_{0<n\in\mathbb{Z}+{1\over 2},\,v>0,\atop n-(p-1)v-\frac{p_j}{n_j}v\leq 0}
(\dim W_v-i_n)(-n+(p-1)v+\tfrac{p_j}{n_j}v)
\\
&\hspace{20pt}+\sum_{0<n\in\mathbb{Z}+{1\over 2},\,v<0,\atop n+(p-1)v+\frac{p_j}{n_j}v\leq 0}
(\dim W_v-i'_n)(-n-(p-1)v-\frac{p_j}{n_j}v)
\\
&=P+(p-1+\tfrac{p_j}{n_j})\emph{\textbf{J}}_H+\varepsilon^2_W\ .
\end{split}
\end{equation}
By \cite[(4.36)-(4.38)]{MR1870666} and \eqref{s42}, we deduce the second line of \eqref{eg68}.
The first line of \eqref{eg68} is obvious.

Consider the $\mathbb{Z}_2$-gradings. By \cite[(4.49)-(4.50)]{MR2016198}
and the discussion following \eqref{uy7}, we get the identities in the first line of \eqref{eg71}.
Observe that $\tau_1$ changes only on
$\bigotimes_{0<n\in\mathbb{Z}+{1\over 2},\,v>0,\atop n-(p-1)v-\frac{p_j}{n_j}v\leq 0}
\Lambda^{i_n}(\overline{W}_{v,n})
\otimes
\bigotimes_{0<n\in\mathbb{Z}+{1\over 2},\,v<0,\atop n+(p-1)v+\frac{p_j}{n_j}v\leq 0}
\Lambda^{i'_n}({W}_{v,n})$. From \eqref{s17}, we get the identity in
the second line of \eqref{eg71}.

The proof of Proposition \ref{s25} is completed.

\subsection{A proof of Theorem \ref{eg31}}\label{sec3.4}

\begin{lem}\label{s30}{\rm (Compare with \cite[Lemmas 4.4 and 4.6]{MR1870666})}
For each connected component $M'$ of $M(n_j)$,
the following functions are independent on the connected components of $F$ in $M'$,
\begin{equation}
\begin{split}\label{s28}
&\varepsilon_i+\varepsilon^1_W+2\varepsilon^2_W\ ,\qquad i=1,2\,,
\\
&d'(p,\beta_{j},N)+\mu_i+\mu_4 \mod(2),\qquad i=1,2,3\,,
\\
&d'(p,\beta_{j-1},N)+\sum_{0<v}\dim N_v+\mu_i+\mu_4 \mod(2),\qquad i=1,2,3\,.
\end{split}
\end{equation}
\end{lem}
{\bf Proof}\hspace{4mm}Recall that $[\tfrac{p_j}{n_j}v]=\tfrac{p_j}{n_j}v-\tfrac{\omega(v)}{n_j}$\ .
By using \eqref{uy3}, we express $\varepsilon^1_W$ and $\varepsilon^2_W$ defined in \eqref{s41}-\eqref{s40}
explicitly as follows,
\begin{align}\label{s31}
\varepsilon^1_W=&\tfrac{1}{2} (p-1+\tfrac{p_j}{n_j})^2 e(W)+\tfrac{1}{8}\dim W(n_j)_{\frac{n_j}{2}}
\notag\\
&+\frac{1}{2}\sum_{0<v<\frac{n_j}{2}}\frac{\omega(v)\omega(-v)}{n_j^2}(\dim W(n_j)_{v}+\dim W(n_j)_{n_j-v})\ ,
\end{align}
\begin{align}\label{s32}
&\varepsilon^2_W=\tfrac{1}{2} (p-1+\tfrac{p_j}{n_j})^2 e(W)-\tfrac{1}{8}\dim W(n_j)_{\frac{n_j}{2}}
\notag\\
&\qquad\quad
-\frac{1}{2}\sum_{0\leq m \leq \frac{p_j-1}{2}}
\sum_{m<\frac{p_j}{n_j}v<m+\frac{1}{2}}\Bigl(\frac{\omega(v)}{n_j}\Bigr)^2
(\dim W(n_j)_{v}+\dim W(n_j)_{n_j-v})
\notag\\
&-\frac{1}{2}\sum_{0\leq m \leq \frac{p_j}{2}}
\sum_{m-\frac{1}{2}<\frac{p_j}{n_j}v<m}\Bigl(\frac{\omega(-v)}{n_j}\Bigr)^2
(\dim W(n_j)_{v}+\dim W(n_j)_{n_j-v})\ .
\end{align}

By using \eqref{uy50}, \eqref{s31}, \eqref{s32} and
the explicit expressions of $\varepsilon_i$
given in \cite[(4.56)-(4.57)]{MR2016198},
we know the functions in the first line of \eqref{s28}
are independent on the connected components of $F$ in $M'$.

Now consider the functions in the rest lines of \eqref{s28}.
By \eqref{s6}, \eqref{uy2}, \eqref{s27}, \eqref{uy7} and \cite[Lemma 4.5]{MR1870666}, we get
\begin{align}
&d'(p,\beta_{j},N)+\mu_i+\mu_4 \equiv
\sum_{0<m\leq \frac{p_j}{2}}
\sum_{0<v<\frac{n_j}{2}\atop m-\frac{1}{2}<\frac{p_i}{n_j}v<m}\dim N(n_j)_{v}
+\frac{1}{2}\dim_\mathbb{R} N(n_j)^\mathbb{R}_{\frac{n_j}{2}}
\notag\\
&\hspace{80pt}+\sum_{v>0}(\dim N_v)\bigl[\tfrac{p_j}{n_j}v+\tfrac{1}{2}\bigr]
+\sum_{v>0}(\dim W_v)\bigl[\tfrac{p_j}{n_j}v+\tfrac{1}{2}\bigr]
\\
&\hspace{20pt}+\sum_{v<0}(\dim W_v)\bigl[-\tfrac{p_j}{n_j}v+\tfrac{1}{2}\bigr]
+o\Bigl(N(n_j)_{\tfrac{n_j}{2}}^{\mathbb{R}}\Bigr)+\Delta(n_j,W) \mod(2)\ .\notag
\end{align}

But by \cite[Lemma 4.5]{MR1870666}, as $w_2(W\oplus TX)_{S^1}=0$, we know that, ${\rm mod}(2)$,
\begin{equation}
\begin{split}
&\sum_{v>0}(\dim N_v)\bigl[\tfrac{p_j}{n_j}v+\tfrac{1}{2}\bigr]
+\sum_{v>0}(\dim W_v)\bigl[\tfrac{p_j}{n_j}v+\tfrac{1}{2}\bigr]
\\
&\hspace{20pt}+\sum_{v<0}(\dim W_v)\bigl[-\tfrac{p_j}{n_j}v+\tfrac{1}{2}\bigr]
+o\Bigl(N(n_j)_{\tfrac{n_j}{2}}^{\mathbb{R}}\Bigr)+\Delta(n_j,W)
\end{split}
\end{equation}
is independent on the connected components of $F$ in $M'$.
Thus, the independence on the connected components of $F$ in $M'$ of
the functions in the second line of \eqref{s28} is proved,
which, combining with \cite[(4.42)]{MR1870666}, implies the
the same independent property of the functions in the third line of \eqref{s28}.

The proof of Lemma \ref{s30} is completed.

\

By \eqref{uy8}-\eqref{s29} and Lemma \ref{eg61}, we know that the Dirac operator
$D^{X(n_j)}\otimes\mathcal{F}(\beta_j)\otimes F_V^i(\beta_j)
\otimes Q_W(\beta_j)\otimes L(\beta_j)_i$ ($i=1,2$) is well-defined on $M(n_j)$.
Observe that \eqref{s43} in Theorem \ref{local} is compatible with the $G_y$ action.
Thus, by using Proposition \ref{s25}, Lemma \ref{s30} and applying Theorem \ref{local} to each connected component of
$M(n_j)$ separately, we deduce that for $i=1,2$, $1\leq j\leq J_0$,
$m\in \tfrac{1}{2}\mathbb{Z}$, $h\in \mathbb{Z}$, $\tau=\tau_{e1}$ or $\tau_{s1}$,
\begin{align}\label{uy11}
&\sum_{\alpha}(-1)^{d'(p,\beta_{j-1},N)+\sum_{v>0}\dim N_v}
\ind_{\tau}\Bigl(D^{Y_\alpha}\otimes(K_W\otimes K^{-1}_X)^{1/2}\otimes\mathcal{F}_{p,j-1}(X)\Bigr.
\notag\\
&\hspace{40pt}
\Bigl.\otimes F_V^i\otimes\ Q(W),m+e(p,\beta_{j-1},N),h\Bigr)
\notag\\
&=\sum_{\beta}(-1)^{d'(p,\beta_{j-1},N)+\sum_{v>0}\dim N_v+\mu}
\ind_{\tau}\Bigl(D^{X(n_j)}\otimes(K_W\otimes K^{-1}_X)^{1/2}\otimes\mathcal{F}(\beta_j)\Bigr.
\notag\\
&\quad \otimes F_V^i(\beta_j)\otimes Q_W(\beta_j)\otimes L(\beta_j)_i,
m+\varepsilon_i+\varepsilon^1_W+2\varepsilon^2_W+(\tfrac{p_j}{n_j}+(p-1))h,h\Big)
\notag\\
&=\sum_{\alpha}(-1)^{d'(p,\beta_{j},N)+\sum_{v>0}\dim N_v}
\ind_{\tau}\Bigl(D^{Y_\alpha}\otimes(K_W\otimes K^{-1}_X)^{1/2}\otimes\mathcal{F}_{p,j}(X)\Bigr.
\notag\\
&\hspace{40pt}
\Bigl.\otimes F_V^i\otimes Q(W),m+e(p,\beta_{j},N),h\Bigr)\ ,
\end{align}
where $\sum_{\beta}$ means the sum over all the connected components of $M(n_j)$. In \eqref{uy11}, if $\tau=\tau_{s1}$, then
$\mu=\mu_3+\mu_4$; if $\tau=\tau_{e1}$, then $\mu=\mu_i+\mu_4$.
Combining \eqref{eg65} with \eqref{uy11}, we get \eqref{eg32}.

The proof of Theorem \ref{eg31} is completed.

\vspace{3mm}\textbf{Acknowledgements}\ \ The authors wish to thank
Professor Weiping Zhang for helpful
discussions.


\end{document}